\documentclass[11pt,a4paper,reqno]{amsart}
\usepackage[T1]{fontenc}
\usepackage[utf8]{inputenc}
\usepackage{lmodern}
\usepackage{amsmath,amssymb,mathtools}
\usepackage[margin=27mm,headheight=14pt]{geometry}
\usepackage{microtype}
\usepackage{enumitem}
\usepackage{needspace}
\usepackage{booktabs,longtable}
\usepackage[colorlinks=true,linkcolor=blue,citecolor=blue,urlcolor=blue]{hyperref}
\hypersetup{pdftitle={A Liouville theorem for the H\'enon parabolic equation},pdfauthor={Yuxia Guo, Wenjun Wang, Tingfeng Yuan},bookmarksnumbered=true}
\newtheorem{theorem}{Theorem}[section]
\newtheorem{proposition}{Proposition}[section]
\newtheorem{lemma}{Lemma}[section]

\theoremstyle{remark}
\newtheorem{remark}{Remark}[section]
\numberwithin{equation}{section}
\allowdisplaybreaks[2]
\newcommand{\R}{\mathbb R}
\newcommand{\Sph}{\mathbb S}
\newcommand{\Ee}{\mathbb E_e}
\newcommand{\Pe}{\mathbb P_e}
\newcommand{\dd}{\,\mathrm d}
\newcommand{\supp}{\operatorname{supp}}
\newcommand{\diver}{\operatorname{div}}

\newcommand{\CC}{\mathcal C}
\newcommand{\KK}{\mathcal K}
\newcommand{\BB}{\mathcal B}

\title[H\'enon parabolic equation]{A Liouville theorem for the H\'enon parabolic equation}
\author{Yuxia Guo}
\author{Wenjun Wang}
\author{Tingfeng Yuan}
\address{Department of Mathematical Sciences, Tsinghua University, Beijing, P.R. China}
\email[Yuxia Guo]{yguo@mail.tsinghua.edu.cn}
\email[Wenjun Wang]{wangwenj25@mails.tsinghua.edu.cn}
\email[Tingfeng Yuan]{ytf22@mails.tsinghua.edu.cn}
\thanks{Corresponding author: Wenjun Wang. E-mail: \href{mailto:wangwenj25@mails.tsinghua.edu.cn}{wangwenj25@mails.tsinghua.edu.cn}.}
\date{}
\subjclass[2020]{Primary 35B53; Secondary 35K58, 35B33, 35B45}
\keywords{H\'enon parabolic equation, Liouville theorem, entire solutions,
Sobolev critical exponent, localized Pohozaev identity}

\begin{document}
\begin{abstract}
We prove that every nonnegative classical entire solution of $u_t-\Delta u=|x|^\alpha u^p$ in $\R^n\times\R$ is identically zero when $n\ge3$, $\alpha\ge0$, and $1<p<\frac{n+2+2\alpha}{n-2}$. The same conclusion holds for weighted parabolic Gross--Pitaevskii systems with nonlinearities of degree $p=2r+1$ and constant coupling coefficients $\beta_{ij}=\beta_{ji}\ge0$ and $\beta_{ii}>0$. Both results include the unweighted case and require no boundedness, radial symmetry, spatial decay, finite-energy, or growth assumptions. The range is optimal for both problems: bounded nonzero nonnegative stationary solutions exist for $p\ge\frac{n+2+2\alpha}{n-2}$.
\end{abstract}
\maketitle

\section{Introduction}\label{1}
We prove a Liouville theorem for the H\'enon parabolic equation
\begin{equation}\label{2}
     u_t-\Delta u=|x|^\alpha u^p \qquad \text{in }\R^n\times\R,\qquad \alpha\ge0,\quad 1 < p < p_S(n, \alpha),
\end{equation}
where
\begin{equation}\label{3}
    p_S(n,\alpha)=
        \begin{cases}
            \frac{n+2+2\alpha}{n-2},&n\ge3,\\[3pt] +\infty,&n=1,2.
        \end{cases}
\end{equation}
When $\alpha=0$, the coefficient in \eqref{2} is understood to be one.

Liouville theorems are closely related to a priori bounds and estimates of singularities and decay; see \cite{21,22,26}. For the unweighted elliptic equation
\begin{equation}\label{4}
    -\Delta u=u^p\qquad\text{in }\R^n,
\end{equation}
Gidas and Spruck~\cite{8} proved that every nonnegative classical solution is zero when $n\ge3$ and $1 < p< p_S(n,0)$. Their proof uses integral estimates based on Bochner identities. Caffarelli, Gidas and Spruck~\cite{3} classified the positive solutions at the critical exponent. Chen and Li~\cite{4} gave proofs of subcritical nonexistence and critical classification using Kelvin inversion and moving planes.

For the heat equation $u_t-\Delta u=u^p$, Fujita \cite{7} proved that positive solutions cannot exist for all forward times if $1<p<1+\frac{2}{n}$. For entire solutions, the expected range was $1<p<p_S(n,0)$. Bidaut-V\'eron \cite{2} extended the elliptic integral
method to space-time and proved the Liouville theorem for
\[
    1<p<p_B(n), \qquad p_B(n)=
    \begin{cases}
        \frac{n(n+2)}{(n-1)^2},&n\ge2,\\[3pt]
        +\infty,&n=1.
    \end{cases}
\]
For $n\ge2$, this left a gap between $p_B(n)$ and $p_S(n,0)$. Merle and Zaag~\cite{13} proved that, for $1<p<p_S(n,0)$, positive ancient solutions satisfying
\[
    u(x,t)\le C(T-t)^{-\frac{1}{p-1}} \qquad(x\in\R^n,\ t<T)
\]
are independent of $x$.

Pol\'a\v cik and Quittner~\cite{19} proved the full subcritical result for bounded nonnegative radial entire solutions by an intersection-comparison argument. Bartsch, Pol\'a\v cik and Quittner~\cite{1} later treated radial sign-changing solutions with a uniformly bounded number of sign changes and derived decay and blow-up estimates. The doubling and rescaling method of Pol\'a\v cik, Quittner and Souplet~\cite{22} reduces the general problem to bounded entire solutions and also gives universal estimates from Liouville theorems. Using energy estimates in backward similarity variables, Quittner~\cite{23} settled the two-dimensional case. This method builds on the similarity variables of Giga and Kohn~\cite{9}. Quittner~\cite{24} then proved the theorem for every $1<p<p_S(n,0)$. His proof improves the rescaled energy estimates in several steps. Local pointwise bounds and the elliptic Liouville theorem are used to obtain these improvements. In the supercritical range, Pol\'a\v cik and Quittner~\cite{20} studied entire and ancient solutions. They proved that every positive bounded radial entire solution is stationary for exponents above the Lepin exponent.

For the weighted elliptic equation
\begin{equation}\label{5}
    -\Delta u=|x|^\alpha u^p\qquad\text{in }\R^n,
\end{equation}
the relevant exponent is $p_S( n, \alpha)$. Phan and Souplet \cite{17} proved the full subcritical result for bounded nonnegative solutions in dimension three using a Pohozaev identity and estimates on spheres. The coefficient in the Pohozaev identity has the required sign precisely when
\[
    \frac{n+\alpha}{p+1}-\frac{n-2}{2}>0 \quad\Longleftrightarrow\quad p<p_S(n,\alpha).
\]
Guo and Wan \cite{10} proved the result for all $n \ge3$. Wang and Ye~\cite{30} obtained Liouville theorems for finite Morse index solutions of H\'enon-type equations in low dimensions and half-spaces; the exponent range in some of their statements was later corrected in a corrigendum. Dai and Qin \cite{5} gave another proof by scaling spheres and extended the method to fractional and higher-order equations. Their comparison with Kelvin transforms gives lower bounds, which are improved by repeated use of the integral equation until a contradiction is reached.

For the weighted parabolic equation \eqref{2}, early work of Filippas and Tertikas \cite{6} treated self-similar solutions. Phan~\cite{15} proved the Liouville theorem for bounded nonnegative entire solutions when $1<p<p_B(n)$, and the full range $1<p<p_S(n,\alpha)$ for nonnegative radial solutions. The proof combines spatial decay, energy dissipation, and the weighted elliptic Liouville theorem. The bounded-solution result covers every $p>1$ when $n=1$. Phan~\cite{16} proved the same conclusion for bounded solutions when $n=2$, using rescaled energy estimates, and obtained partial results for $n\ge3$. For $n=1,2$, doubling and rescaling also remove boundedness: the limits solve either a translated weighted equation or an unweighted equation, and the bounded Liouville theorems exclude both. For $n\ge3$, the question is to prove the result throughout $1<p<p_S(n,\alpha)$ without radial symmetry, as conjectured in \cite[Remark~1.1(c)]{15}.

For systems, Serrin and Zou~\cite{28} established early nonexistence results for positive solutions of the Lane--Emden system. Reichel and Zou~\cite{27} used moving spheres to prove nonexistence results for cooperative semilinear elliptic systems. Souplet~\cite{29} proved the Lane--Emden conjecture in dimension four and obtained further nonexistence ranges in higher dimensions. For weighted Lane--Emden systems, Hu~\cite{11} proved nonexistence results for positive semistable solutions. Phan and Souplet~\cite{18} proved a Liouville theorem for the unweighted version of \eqref{10} when $p=2r+1<p_B(n)$, with the same assumptions on the matrix as in Theorem~\ref{9} below. Their integral method includes the cubic case $r=1$ in dimension three. Quittner~\cite{24, 25} obtained parabolic Liouville theorems for homogeneous systems with gradient structure from the corresponding elliptic theorems, for $1<p<p_S(n,0)$.

For the elliptic H\'enon--Lane--Emden system
\[
    \left\{
        \begin{aligned}
        -\Delta u&=|x|^A v^p,\\
        -\Delta v&=|x|^B u^q,
        \end{aligned}
    \right.\qquad\text{in }\R^n\setminus\{0\},
\]
OpenAI \cite[Theorem~1.1]{14} proved that this system has no solutions $u,v\in C^2(\R^n\setminus\{0\})\cap C(\R^n)$ with $u,v>0$ on $\R^n$ when
\[
    n\ge2,\qquad p,q>0,\qquad A,B\in\R,\qquad \frac{n+A}{p+1}+\frac{n+B}{q+1}>n-2.
\]
Taking $A=B=0$ gives a proof of the Lane--Emden conjecture. The comparison on intervals, angular localization, and integration over source levels in that proof suggested the approach used here.

The weighted parabolic problem has two difficulties. First, the weighted elliptic theorem does not apply directly at a fixed time, since
\[
    -\Delta u(\cdot,t)=|x|^\alpha u(\cdot,t)^p-u_t(\cdot,t),
\]
and the last term has no fixed sign.

Second, Quittner's proof \cite[Lemma~6]{24} obtains local pointwise bounds from integral bounds on the time derivative. Doubling and rescaling give a positive stationary limit, which is excluded by the unweighted elliptic Liouville theorem. For the weighted equation, a normalized rescaling near $x_j\ne0$ at a scale $r_j$ with $r_j/|x_j| \to 0 $ has coefficient
\begin{equation}\label{6}
    \left(\frac{|x_j+r_jy|}{|x_j|}\right)^\alpha \to 1 \quad\text{locally uniformly in }y, \ \ \text{as} \ \  j \to +\infty.
\end{equation}
A stationary limit would therefore solve $-\Delta v=v^p$.
This equation has bounded positive solutions when
\begin{equation}\label{7}
    p_S(n,0)\le p<p_S(n,\alpha).
\end{equation}
Thus this argument for pointwise estimates does not apply to the full weighted subcritical range.

Our first result covers the full weighted subcritical range.
\begin{theorem}\label{8}
    Let $n\ge3$, $\alpha\ge0$, and $1<p<p_S(n,\alpha)$. Then every nonnegative solution $ u \in C^{2,1}(\R^n \times \R)$ of \eqref{2} satisfies $u \equiv 0$.
\end{theorem}

The same method applies to weighted parabolic Gross--Pitaevskii
and related systems.
\begin{theorem}\label{9}
    Let $n\ge3$, $\alpha\ge0$, $m\ge2$, and $r>0$. Let $(\beta_{ij})$ be a constant real symmetric $m\times m$ matrix such that
    \[
        \beta_{ij}\ge0\quad(1\le i,j\le m),\qquad \beta_{ii}>0\quad(1\le i\le m),\qquad 1<p:=2r+1<p_S(n,\alpha).
    \]
    If $\mathbf u=(u_1,\ldots,u_m)\in C^{2,1}(\R^n\times\R;[0,\infty)^m)$ solves
    \begin{equation}\label{10}
        (u_i)_t-\Delta u_i =|x|^\alpha u_i^r\sum_{j=1}^m\beta_{ij}u_j^{r+1} \quad\text{in }\R^n\times\R,\qquad i=1,\ldots,m,
    \end{equation}
    then $\mathbf u\equiv0$.
\end{theorem}
The proof uses the same method as the proof of Theorem~\ref{8}. A sketch is given in Section~\ref{123}.

Neither theorem requires boundedness, radial symmetry, spatial decay, finite energy, or a growth bound at infinity. For $\alpha=0$, Theorem \ref{8} includes Quittner's scalar result in dimensions $n\ge3$. The cubic Gross--Pitaevskii case is $r=1$, for which the condition
in Theorem \ref{9} becomes $n<4+\alpha$.

\begin{remark}\label{11}
    The range is optimal for both theorems. At $p=p_S(n,\alpha)$,
    \begin{equation}\label{12}
        U(x)=\frac{\bigl[(n-2)(n+\alpha)\bigr]^{\frac{n-2}{2(\alpha+2)}}} {(1+|x|^{\alpha+2})^{\frac{n-2}{\alpha+2}}}
    \end{equation}
    is a bounded positive stationary solution of \eqref{2}. Such radial solutions also exist for every $p>p_S(n,\alpha)$; see \cite[Remark~1.1(a)]{15}. If $U$ is any of these solutions, then
    \[
        \mathbf u(x,t)= \bigl(\beta_{11}^{-\frac{1}{p-1}}U(x),0,\ldots,0\bigr)
    \]
    is a bounded nonzero nonnegative solution of \eqref{10}.
\end{remark}

In our proof, we use a resolvent potential to separate the time derivative from the positive source. A local pressure estimate, Pohozaev identities, and an iteration on nested cylinders give
\[
    R^\delta\iint_{Q_R}u_t^2 +R^{\delta-2}\iint_{Q_R}  \bigl(|\nabla u|^2+|x|^\alpha u^{p+1}\bigr) \le C(n,p,\alpha),
\]
where
\[
    Q_R=B_R\times(-R^2,R^2),\qquad \delta=\frac{2(\alpha+2)}{p-1}+2-n>0.
\]
Letting $R\to\infty$ first gives $u_t=0$. The same estimate then gives $u=0$.

The paper is organized as follows. Section \ref{13} establishes the mass and interpolation estimates and derives a localized Pohozaev identity. Section \ref{41} develops the resolvent and interval comparison methods and proves the local pressure estimate. Section \ref{123} derives a uniform local energy bound and completes the proof of Theorem \ref{8}. Throughout the paper, $C,c$ denote positive constants that may vary from line to line. 

\section{Preliminary estimates and a localized Pohozaev identity}\label{13}
Let $n\ge3$, $\alpha\ge0$, and $p>1$. We consider nonnegative $C^{2,1}$ solutions of
\begin{equation}\label{14}
    u_t-\Delta u = |x|^\alpha u^p\qquad\text{in } \R^n \times \R.
\end{equation}
Set
\begin{equation}\label{15}
    \kappa=\frac{\alpha+2}{p-1},\qquad \delta=2\kappa+2-n.
\end{equation}

We use Kaplan's weighted-mean argument \cite{12} on an annulus to obtain the following lemma.

\begin{lemma}\label{16}
Suppose $u$ is a nonnegative classical solution of \eqref{14}, and let $ \kappa $ be as in \eqref{15}. Then there is a constant $ C_0 = C_0(n,p,\alpha) > 0 $ such that
\begin{equation}\label{17}
    \int_{B_R} u(x,t) \dd x \le C_0 R^{n-\kappa}\qquad(R>0,\ t\in\R).
\end{equation}
\end{lemma}

\begin{proof}
    Let $A_1=\{x\in\R^n:1<|x|<2\}$, and let $\phi_1>0$ be the first Dirichlet eigenfunction on $A_1$, normalized by $\int_{A_1}\phi_1=1$:
    \[
        \left\{
        \begin{aligned}
        -\Delta\phi_1&=\lambda_1\phi_1 && \text{in } A_1,\\
        \phi_1 & = 0 && \text{on } \partial A_1,
        \end{aligned} \right. \qquad \lambda_1>0.
    \]
    For the outward unit normal $\nu$, the Hopf lemma gives $ \partial_\nu \phi_1 < 0 $ on $ \partial A_1 $. Set
    \[
        A'_1= \left\{ \frac{5}{4}<|x|<\frac{7}{4} \right\} , \qquad \min_{\overline{A'_1}}\phi_1\ge c>0.
    \]
    For $R>0$, set
    \[
        A_R=\{R<|x|<2R\},\qquad \phi_R(x)=R^{-n} \phi_1 \left(\frac{x}{R} \right),\qquad F_R(t)=\int_{A_R}u(x,t) \phi_R(x)\dd x.
    \]
    Since $ \int_{A_R} \phi_R = 1 $, $|x|^\alpha \ge R^\alpha$ on $A_R$, and $ \partial_\nu \phi_R < 0 $, Green's formula and Jensen's inequality give
    \begin{equation}\label{18}
    \begin{aligned}
        F_R' & = -\lambda_1 R^{-2} F_R -\int_{\partial A_R} u \partial_\nu \phi_R + \int_{A_R}|x|^\alpha u^p\phi_R \\
        & \ge - \lambda_1 R^{-2} F_R + R^\alpha F_R^p.
    \end{aligned}
    \end{equation}
    We conclude from \eqref{18} that
    \begin{equation}\label{19}
        F_R(t) \le \lambda_1^{\frac{1}{p-1}} R^{-\kappa} \qquad (t \in \R).
    \end{equation}
    Suppose, for contradiction, that some $t_0 \in \R$ satisfies 
    \begin{align}\label{20}
        F_R(t_0) > \lambda_1^{\frac{1}{p-1}}R^{-\kappa}.
    \end{align}
    By \eqref{18}, $F_R$ remains above this threshold and $F_R(t)\ge F_R(t_0)$ for $t\ge t_0$. Multiplication by $(1-p)F_R^{-p}<0$ gives
    \[
        \frac{\mathrm d}{\mathrm dt}F_R^{1-p} \le(p-1)\lambda_1R^{-2}F_R^{1-p}-(p-1)R^\alpha.
    \]
    Integration yields
    \[
        F_R(t)^{1-p}\le\frac{R^{\alpha+2}}{\lambda_1} + \left(F_R(t_0)^{1-p}-\frac{R^{\alpha+2}}{\lambda_1}\right) e^{(p-1)\lambda_1R^{-2}(t-t_0)}.
    \]
    By \eqref{20}, the right-hand side reaches zero in finite time, contradicting $F_R(t)^{1-p}>0$. Hence, \eqref{19} holds.

    Set $ A'_R = R A'_1 = \{x\in\R^n:\frac{5R}{4}<|x|< \frac{7R}{4} \}$. Since $\phi_R\ge cR^{-n}$ on $A'_R$, we obtain from \eqref{19} that
    \begin{equation}\label{21}
        \int_{A_R'} u(x,t)\dd x\le CR^{n-\kappa}.
    \end{equation}
    Finally, $u_t-\Delta u\ge0$. We claim that
    \begin{equation}\label{22}
        u(x,t+s)\ge\int_{B_R}\Gamma(x-y,s)u(y,t)\dd y\qquad(s>0),
    \end{equation}
    where $\Gamma(z,s)=(4\pi s)^{-\frac{n}{2}}e^{-\frac{|z|^2}{4s}}$ is the heat kernel. Choose $\eta_j\in C_c^\infty(B_R)$, extended by zero outside $B_R$, such that
    \[
        0\le\eta_j(x)\le\eta_{j+1}(x)\le1,\qquad \lim_{j\to\infty}\eta_j(x)=\mathbf1_{B_R}(x) \quad(x\in\R^n).
    \]
    Here $\mathbf1_A$ denotes the characteristic function of $A$. The difference
    \[
        [\Gamma(\cdot,s)*(\eta_j u(\cdot,t))](x)-u(x,t+s) \le\|\eta_j u(\cdot,t)\|_\infty
    \]
    is a subsolution of the heat equation that is bounded above and has nonpositive initial values. The maximum principle \cite[Proposition~52.4]{26} and the monotone convergence theorem give \eqref{22}. For $x\in A'_R$ and $y\in B_R$, $|x-y|\le\frac{11R}{4}$ and $\Gamma(x-y,R^2)\ge c_nR^{-n}$. Thus
    \[
        \int_{B_R} u(y,t) \dd y \le C_n \int_{A'_R} u(x,t+R^2)\dd x \le CR^{n-\kappa},
    \]
    by \eqref{21}. This proves \eqref{17}.
\end{proof}

For $a>0$, let
\begin{equation}\label{23}
    v_a(y,s)=a^\kappa u(ay,a^2s).
\end{equation}
Since $\kappa(p-1)=\alpha+2$, the rescaled function satisfies
\[
    \partial_sv_a-\Delta_yv_a=|y|^\alpha v_a^p.
\]
For every $R>0$, Lemma~\ref{16} gives
\begin{equation}\label{24}
\begin{aligned}
    \sup_s\int_{B_R}v_a(y,s)\dd y = a^{\kappa-n} \sup_s\int_{B_{aR}}u(x,a^2s)\dd x \le C_0 R^{n-\kappa}.
\end{aligned}
\end{equation}

Fix $a>0$ and write $v=v_a$, with variables
denoted by $(x,t)$. Then
\begin{equation}\label{25}
    v_t-\Delta v=|x|^\alpha v^p.
\end{equation}
For every $\rho>0$, \eqref{24} gives
\[
    \sup_{t\in\R}\int_{B_\rho}v(x,t)\dd x\le C_0\rho^{n-\kappa}.
\]
Define
\begin{equation}\label{26}
    E_\rho(t) = \int_{B_\rho}( |\nabla v|^2 + |x|^\alpha v^{p+1} )\dd x, \qquad D_\rho(t) = \int_{B_\rho} v_t^2 \dd x.
\end{equation}
Since $v\in C^{2,1}$, both functions are continuous on $\R$ and finite at every time. We fix the time $ t \in \R $ and suppress it in the spatial estimates.

\begin{lemma}\label{27}
    At each fixed $ t \in \R $,
    \begin{equation}\label{28}
        \|v\|_{L^2(B_\rho)}^2 \le C_n C_0^{\frac{4}{n+2}} \rho^{\frac{4(n-\kappa)}{n+2}} \left( \int_{B_\rho} |\nabla v|^2 \right)^{\frac{n}{n+2}} + C_n C_0^2 \rho^{n-2\kappa}.
    \end{equation}
    Consequently, for $ 0 < \varepsilon \le 1 $,
    \begin{equation}\label{29}
        \|v\|_{L^2(B_\rho)}^2 \le \varepsilon \int_{B_\rho} |\nabla v|^2 + C_n C_0^2 \bigl( \varepsilon^{-\frac{n}{2}} \rho^{2(n-\kappa)} + \rho^{n-2\kappa}\bigr).
    \end{equation}
    Set
    \begin{equation}\label{30}
        N_\rho = 1 + \|v\|_{L^2(B_\rho)}^2 + \|v\|_{L^2(B_\rho)}(\sqrt{E_\rho} + \sqrt{D_\rho}) + \sqrt{D_\rho}.
    \end{equation}
    For the integer $ m_2 = 2n+5 $, we have
    \begin{align}
        N_\rho & \le C_n( 2 +C_0\rho^{n-\kappa}+\rho^{-1} )^{n+2} \bigl[ 1+(E_\rho+D_\rho)^{\frac{n+1}{n+2}} \bigr], \label{31}\\
        L N_\rho & \le \varepsilon( E_\rho + D_\rho) + C_n( 2 + L + C_0 \rho^{n-\kappa} + \rho^{-1} + \varepsilon^{-1})^{m_2}, \quad L \ge 0. \label{32}
    \end{align}
\end{lemma}

\begin{proof}
    H\"older interpolation and the Sobolev--Poincar\'e inequality on $B_1$ give
    \[
        \begin{aligned}
            \|v\|_{L^2(B_1)} & \le \|v\|_{ L^1(B_1) }^{\frac{2}{n+2}}  \|v\|_{L^{\frac{2n}{n-2}} (B_1)}^{\frac{n}{n+2}}, \\
            \|v\|_{L^{\frac{2n}{n-2}} (B_1)} & \le C_n \bigl( \| \nabla v \|_{L^2 (B_1)} + \|v\|_{L^1(B_1)}\bigr).
        \end{aligned}
    \]
    Combining these inequalities and scaling $x=\rho y$ yields
    \begin{equation}\label{33}
        \|v\|_{L^2(B_\rho)}^2 \le C_n \left[\left( \int_{B_\rho} v \right)^{\frac{4}{n+2}} \left(\int_{B_\rho}| \nabla v|^2\right)^{\frac{n}{n+2}} + \rho^{-n} \left( \int_{B_\rho} v \right)^2\right].
    \end{equation}
    Substitute $\int_{B_\rho} v \le C_0 \rho^{n-\kappa}$ to get \eqref{28}. Young's inequality with exponents $\frac{n+2}{n}$ and $\frac{n+2}{2}$ gives \eqref{29}. Taking square roots in \eqref{28} gives
    \[
        \|v\|_{L^2(B_\rho)} \le C_nC_0^{\frac{2}{n+2}} \rho^{\frac{ 2(n-\kappa) }{n+2}} E_\rho^{\frac{n}{2(n+2)}} + C_n C_0 \rho^{\frac{n}{2}-\kappa}.
    \]
    Substitution in \eqref{30} gives
    \[
        N_\rho \le C_n( 2 + C_0 \rho^{n-\kappa} + \rho^{-1} )^{n+2} [1+(E_\rho+D_\rho)^{\frac{n+1}{n+2}}].
    \]
    Applying Young's inequality separately to the terms in $L N_\rho$ gives
    \[
        \begin{aligned}
            L N_\rho & \le \varepsilon(E_\rho+D_\rho) + C_n( L + \varepsilon^{-1} L^2)(1+C_0^2\rho^{n-2\kappa}) \\
            & \quad + C_n\left[ \varepsilon^{-\frac{n}{2}} L^{\frac{n+2}{2}} + \varepsilon^{-(n+1)} L^{n+2}\right] C_0^2\rho^{2(n-\kappa)}\\
            & \le \varepsilon(E_\rho+D_\rho) + C_n( 2+ L+ C_0\rho^{n-\kappa} + \rho^{-1} + \varepsilon^{-1} )^{2n+5}.
        \end{aligned}
    \]
    This proves \eqref{31} and \eqref{32}.
\end{proof}

\begin{remark}\label{34}
    The proof uses only the mass bound and spatial interpolation. Hence the same estimates hold for local classical solutions under the mass bound assumed in Proposition~\ref{75}.
\end{remark}

For a radial function $\varphi\in C_c^3(B_\rho)$, write $\varphi(x)=\varphi(r)$, where $r=|x|$. Set
\[
    w(r)=-r\varphi'(r),\qquad e_r(x)=\frac{x}{|x|},\qquad v_r = e_r \cdot \nabla v \quad( x \ne 0 ).
\]
The following identity relates the nonlinear term to a radial pressure.
\begin{lemma}\label{35}
    For a solution of \eqref{25}, at each fixed time,
    \begin{equation}\label{36}
    \begin{aligned}
        &\int\varphi v_t\left( x \cdot \nabla v + \frac{n-2}{2} v \right) + \left( \frac{ n+\alpha }{p+1} - \frac{n-2}{2} \right) \int \varphi |x|^\alpha v^{p+1}\\
        & \quad = \int w \left( v_r^2 - \frac{p-1}{2(p+1)} |x|^\alpha v^{p+1}\right) + \frac{1}{2}\int w v v_t + \frac{1}{4} \int \bigl((n-2)\Delta\varphi-\Delta w \bigr)v^2.
    \end{aligned}
    \end{equation}
Moreover,
    \begin{equation}\label{37}
    \begin{aligned}
        & \left| \frac{1}{2} \int w v v_t\right| + \left| \frac{1}{4} \int \bigl( (n-2) \Delta \varphi -\Delta w \bigr) v^2 \right|\\
        \quad \le & \frac{1}{2} \|w\|_\infty \|v\|_{L^2(B_\rho)} \sqrt{D_\rho} + \frac{1}{4} \|(n-2) \Delta \varphi - \Delta w \|_\infty \|v\|_{L^2(B_\rho)}^2 \\
        \quad \le & C_{n,\varphi} \bigl( \|v \|_{L^2(B_\rho)} \sqrt{D_\rho} + \| v \|_{L^2(B_\rho)}^2 \bigr).
    \end{aligned}
    \end{equation}
\end{lemma}

\begin{proof}
    Fix a time $t$. Since $\nabla \varphi = \frac{\varphi'(r)x}{r} $ and $ r \varphi'(r) = -w $, integration by parts gives
    \begin{align}\label{38}
        \int \varphi\Delta v\,x\cdot\nabla v &=\frac{n-2}{2} \int \varphi| \nabla v|^2 + \int w\left(v_r^2-\frac{1}{2}|\nabla v|^2\right).
    \end{align}
    
    For the nonlinear term, use $v^p x \cdot \nabla v = \frac{1}{p+1} x \cdot \nabla(v^{p+1})$ and $\operatorname{div} (|x|^\alpha x)=(n+\alpha) |x|^\alpha$ to obtain
    \begin{equation}\label{39}
        \int\varphi |x|^\alpha v^p x\cdot\nabla v = - \frac{n+\alpha}{p+1}\int\varphi|x|^\alpha v^{p+1} + \frac{1}{p+1}\int w|x|^\alpha v^{p+1}.
    \end{equation}
    Testing \eqref{25} with $\varphi v$ also gives
    \begin{equation}\label{40}
        \int\varphi v v_t = \int\varphi(|x|^\alpha v^{p+1}-|\nabla v|^2) + \frac{1}{2} \int(\Delta\varphi)v^2.
    \end{equation}
    Add $\frac{n-2}{2}$ times \eqref{40} to the sum of \eqref{38} and \eqref{39}. The gradient terms with weight $\varphi$ cancel. In the remaining integral, write
    \[
        v_r^2-\frac{1}{2}|\nabla v|^2 + \frac{|x|^\alpha v^{p+1}}{p+1} = v_r^2 - \frac{p-1}{2(p+1)} |x|^\alpha v^{p+1} + \frac{1}{2} (|x|^\alpha v^{p+1}-|\nabla v|^2).
    \]
    Applying \eqref{40} with $w$ in place of $\varphi$ gives
    \[
        \int w(|x|^\alpha v^{p+1}-|\nabla v|^2) = \int wvv_t - \frac{1}{2} \int(\Delta w)v^2.
    \]
    Substitution proves \eqref{36}. The bound \eqref{37} follows from the Cauchy--Schwarz inequality. \end{proof}

The identity holds for every $p>1$, and
\[
    \frac{n+\alpha}{p+1}-\frac{n-2}{2} = \frac{p-1}{2(p+1)}\delta>0 \qquad\text{if }p<p_S(n,\alpha).
\]

\section{Resolvent methods and a local pressure estimate}\label{41}
   In this section, we prove the local pressure estimate used in Section~\ref{123}. We begin with resolvent estimates and an interval comparison.

   Fix $R_*>0$, $0<r_0<r_1$, and $0<h\le1$ such that $\rho=r_1+6h\le R_*$. Let $w_0\in C_c^2(\{r_0<|x|<r_1\})$ be nonnegative and radial. We work at a fixed time $t$ at which $v$ is a nonnegative classical solution of \eqref{25} near $\overline{B_\rho}$ and satisfies
    \[
        \int_{B_\rho}v(x,t)\dd x\le C_0\rho^{n-\kappa}.
    \]
    We suppress the time variable in the estimates below. For $e\in\Sph^{n-1}$, write
    \[
        \partial_e=e\cdot\nabla,\qquad v_e = \partial_e v,\qquad v_{ee}=\partial_e^2v, \qquad \|v\|_2  = \| v(\cdot,t) \|_{L^2(B_\rho)}.
    \]
    We use the same convention for other directional derivatives. Norms of extended functions are over $\R^n$. The quantities $E_\rho,D_\rho,N_\rho$ are defined in \eqref{26} and \eqref{30}. Put
    \[
        Q_0 = 2 + C_0 \rho^{n-\kappa} + h^{-1} + r_0^{-1} + \rho^{-1} + \|w_0\|_{C^2},
    \]
    where $C_0$ is the constant from Lemma~\ref{16}.

    \subsection{Resolvent and interval estimates.} For $\lambda\ge1$, let $R_\lambda=(-\Delta+\lambda)^{-1}$ be the resolvent on $\R^n$. Thus $R_\lambda q=K_\lambda*q$, where
    \begin{equation}\label{42}
        K_\lambda(z)= \int_0^\infty e^{-\lambda a} \Gamma(|z|,a) \dd a = \int_0^\infty e^{-\lambda a}(4\pi a)^{-\frac{n}{2}} e^{-\frac{ |z|^2}{4a}} \dd a .
    \end{equation}
    The kernel is positive, radial, and decreasing in $|z|>0$. We also write $K_\lambda(z)=K_\lambda(r)$, where $r=|z|$. Gaussian integration gives
\begin{equation}\label{43}
    \|K_\lambda\|_1=\lambda^{-1},\qquad \widehat K_\lambda(\xi) = \frac{(2\pi)^{-\frac{n}{2}}}{\lambda+|\xi|^2}.
\end{equation}

Thus $R_\lambda$ is positive and self-adjoint on $L^2$ and commutes with derivatives. For a multi-index $\gamma$, Gaussian differentiation and scaling give
\begin{align}
    |\partial^\gamma\Gamma(x,t)| & \le C_\gamma t^{-\frac{n+|\gamma|}{2}} e^{-\frac{|x|^2}{8t}}, \label{44}\\
    \|\partial^\gamma\Gamma(\cdot,t)\|_1 & \le C_\gamma t^{-\frac{|\gamma|}{2}}, \label{45}
\end{align}
where the second inequality follows from the substitution $x=\sqrt t\,z$.

For $0\le j\le4$ and $z\ne0$, these estimates give
\begin{equation}\label{46}
    |D^jK_\lambda(z)| \le C_j \int_0^\infty a^{-\frac{n+j}{2}} e^{-\frac{|z|^2}{8a}} \dd a \le C_{n,j} |z|^{2-n-j}.
\end{equation}
The integral is finite since $n\ge3$. In particular, the first distributional derivatives are integrable:
\[
    \|\nabla K_\lambda\|_1 \le C\int_0^\infty e^{-\lambda t} t^{-\frac{1}{2}} \dd t \le C\lambda^{-\frac{1}{2}}.
\]

For $r>0$, differentiation gives
\begin{equation}\label{47}
    K_\lambda'(r)=-\int_0^\infty e^{-\lambda t}\frac{r}{2t} \Gamma(r,t)\dd t<0.
\end{equation}
The bounds above justify this differentiation. The substitution $s=\lambda t$ yields
\begin{equation}\label{48}
    K_\lambda(x)=\lambda^{\frac{n-2}{2}} K_1(\sqrt\lambda\,x).
\end{equation}
For $n\ge3$, the same Gaussian integral with $\lambda=0$ gives the Newton kernel
\begin{equation}\label{49}
    K_0(x)=\frac{|x|^{2-n}}{(n-2)|\Sph^{n-1}|}.
\end{equation}
Positivity of the integral gives $K_\lambda\le K_0$.

For $0<r=|x|\le\lambda^{-\frac{1}{2}}$, restrict the integral defining $K_\lambda$ to $r^2\le t\le2r^2$ and put $t=r^2s$. Then
\begin{equation}\label{50}
    K_\lambda(x)\ge |x|^{2-n} \int_1^2e^{-2}(4\pi s)^{-\frac{n}{2}}e^{-\frac{1}{4s}}\dd s \ge c(n)|x|^{2-n},
\end{equation}
where $c(n)>0$ depends only on the dimension.

For $r=|x|\ge1$,
\[
    t+\frac{r^2}{8t}\ge\frac{t}{2}+\frac{1}{16t}+\frac{r}{4}.
\]
Thus \eqref{44} implies $|\partial^\gamma K_\lambda(x)|\le C_\gamma e^{-\frac{|x|}{4}}$ for $|x|\ge1$. For smooth compactly supported sources, the resolvent and its derivatives therefore decay fast enough that boundary terms from integration by parts on expanding balls tend to zero.

Choose $\eta\in C_c^\infty(B_{r_1+5h})$, $0\le\eta\le1$, equal to one on $B_{r_1+4h}$, with $\|D^j\eta\|_\infty\le C_jh^{-j}$ for $1\le j\le4$. Set
\begin{equation}\label{51}
    g=\eta(|x|^\alpha v^p+\lambda v).
\end{equation}
Extend localized functions by zero outside $B_\rho$. The equation gives
\begin{equation}\label{52}
    g=(-\Delta+\lambda)(\eta v)+\eta v_t + 2 \operatorname{div} (v\nabla\eta) - v\Delta\eta.
\end{equation}
Consequently,
\begin{equation}\label{53}
\begin{gathered}
    W=R_\lambda g=\eta v+H+B,\qquad H=R_\lambda(\eta v_t),\\
    B=R_\lambda\bigl(2\operatorname{div}(v\nabla\eta)-v\Delta\eta\bigr), \qquad \CC=H+B.
\end{gathered}
\end{equation}
Here $W\ge0$ and $W=v+\CC$ on $B_{r_1+2h}$. For $0 < h \le 1 $, we claim that
\begin{equation}\label{54}
\begin{aligned}
    \|B\|_2 &\le C(h^{-1}\lambda^{-\frac{1}{2}}+h^{-2}\lambda^{-1})\|v\|_2,\\
    \|\nabla B\|_2 &\le C(h^{-1}+h^{-2}\lambda^{-\frac{1}{2}})\|v\|_2,\\
    \|B\|_{C^2(B_{r_1+2h})} &\le Ch^{-n-2}C_0\rho^{n-\kappa}.
\end{aligned}
\end{equation}
First, we have $\|\eta v_t\|_2^2\le D_\rho$ and
\[
    \widehat H(\xi)=\frac{\widehat{\eta v_t}(\xi)}{\lambda+|\xi|^2}.
\]
Plancherel's identity gives
\begin{align*}
    \|H\|_2^2 &=\int_{\R^n}\frac{|\widehat{\eta v_t}(\xi)|^2}{(\lambda+|\xi|^2)^2}\dd\xi \le \lambda^{-2}D_\rho,\\
    \|\nabla H\|_2^2 & = \int_{\R^n} \frac{|\xi|^2| \widehat{\eta v_t}(\xi)|^2}{(\lambda+|\xi|^2)^2} \dd\xi \le\frac{D_\rho}{4\lambda},\\
    \|D^2H\|_2^2 &=\int_{\R^n}\frac{|\xi|^4|\widehat{\eta v_t}(\xi)|^2}{(\lambda+|\xi|^2)^2}\dd\xi \le D_\rho.
\end{align*}
Indeed, for $r\ge0$,
\[
    \frac{1}{\lambda+r^2}\le\lambda^{-1},\qquad \frac{r}{\lambda+r^2}\le\frac{1}{2\sqrt\lambda},\qquad \frac{r^2}{\lambda+r^2}\le1.
\]
Thus
\begin{equation}\label{55}
    \|H\|_2\le\lambda^{-1}\sqrt{D_\rho},\qquad \|\nabla H\|_2\le\frac{\sqrt{D_\rho}}{2\sqrt\lambda},\qquad \|D^2H\|_2\le\sqrt{D_\rho}.
\end{equation}
For a vector field $F\in L^2(\R^n;\R^n)$, the same multipliers give
\[
    \|R_\lambda\diver F\|_2\le\frac{\|F\|_2}{2\sqrt\lambda}, \qquad \|\nabla R_\lambda\diver F\|_2\le\|F\|_2.
\]
Applying these inequalities to $F=2v\nabla\eta$ and using
$\|v\Delta\eta\|_2\le Ch^{-2}\|v\|_2$ yields the $L^2$ bounds for $B$ and $\nabla B$.

To obtain the local $C^2$ bound, integration by parts in the convolution gives
\[
    B(x)=2\int\nabla_xK_\lambda(x-y)\cdot v(y)\nabla\eta(y)\dd y -\int K_\lambda(x-y)v(y)\Delta\eta(y)\dd y.
\]
By the choice of $\eta$,
\[
    \supp\nabla\eta\cup\supp\Delta\eta \subset B_{r_1+5h}\setminus B_{r_1+4h}.
\]
Hence, if $x\in B_{r_1+2h}$ and $y\in\supp\nabla\eta\cup\supp\Delta\eta$, then $|x-y|\ge2h$. Differentiation under the integrals is therefore legitimate, and \eqref{46} gives, for $j=0,1,2$,
\[
    |D^jB(x)| \le Ch^{-n-j}\int_{B_\rho}v(y)\dd y \le Ch^{-n-j} C_0\rho^{n-\kappa}.
\]
These estimates prove \eqref{54}. In particular, these bounds are uniform for $\lambda\ge1$. Moreover,
\[
    \|\eta v\|_2\le\|v\|_2,\qquad \|\nabla(\eta v)\|_2\le\sqrt{E_\rho}+Ch^{-1}\|v\|_2.
\]
The definition of $W$ gives, for $\lambda\ge1$,
\begin{equation}\label{56}
\begin{aligned}
    \|W\|_2 & \le C(1+h^{-2}) \|v\|_2+ \lambda^{-1} \sqrt{D_\rho},\\
    \|\nabla W\|_2 & \le \sqrt{E_\rho} + C(1+h^{-2})\|v\|_2 + \frac{\sqrt{D_\rho}}{2\sqrt\lambda}.
\end{aligned}
\end{equation}
Testing \eqref{25} with $\eta$ gives
\begin{equation}\label{57}
\begin{aligned}
    \int g &=\lambda\int\eta v+\int\eta v_t-\int v\Delta\eta\\
    & \le (\lambda+Ch^{-2}) C_0 \rho^{n-\kappa} +C_n\rho^{\frac{n}{2}}\sqrt{D_\rho}.
 \end{aligned}
\end{equation}
Positivity and self-adjointness of $R_\lambda$ imply
\begin{equation}\label{58}
    \int_{B_{r_1+2h}}g|H| \le\int gR_\lambda|\eta v_t| =\int|\eta v_t|W \le C(1+h^{-2}) \|v\|_2 \sqrt{D_\rho} +\lambda^{-1}D_\rho.
\end{equation}
Combining this with \eqref{54} and \eqref{57}, we obtain
\begin{equation}\label{59}
 \begin{aligned}
    \int_{B_{r_1+2h}}g|\CC| & \le C(1+h^{-2}) \|v\|_2\sqrt{D_\rho} +\lambda^{-1}D_\rho\\
     &\quad+Ch^{-n}C_0\rho^{n-\kappa} \left[(\lambda+ h^{-2}) C_0 \rho^{n-\kappa} + \rho^{\frac{n}{2}}\sqrt{D_\rho}\right].
\end{aligned}
\end{equation}

Note that the estimates \eqref{55}--\eqref{59} require no sign condition on $v_t$. We will use them repeatedly below.

The following interval comparison is adapted from \cite{14}. We give its proof in the form needed here.

Fix $e\in\Sph^{n-1}$ and let $e^\perp=\{z\in\R^n:z\cdot e=0\}$. For $z,z'\in e^\perp$ with $z\ne z'$, write
$x=z+se$, $y=z'+s'e$, and
\[
    k(s-s')=K_\lambda(z-z'+(s-s')e).
\]
Since $z-z'\perp e$,
\[
    k(\xi)=K_\lambda(\sqrt{|z-z'|^2+\xi^2})
\]
is smooth, positive, and even, with
\[
    k'(\xi)=K_\lambda'(\sqrt{|z-z'|^2+\xi^2}) \frac{\xi}{\sqrt{|z-z'|^2+\xi^2}}\le0 \quad(\xi>0).
\]
For intervals $I=(a,b)$ and $J=(a',b')$ with $a<b$ and $a'<b'$, let
\[
    k_I(s')=\frac{k(a-s')+k(b-s')}{2},\qquad k_J(s)=\frac{k(s-a')+k(s-b')}{2}.
\]

\begin{lemma}\label{60}
For $c_I, c_J \ge0$,
\begin{align}
    \int_I\int_J(c_I+c_J)k(s-s')\dd s'\dd s \le \int_I \int_J \bigl[& c_Ik_I(s')+c_Jk_J(s) -\min(c_I,c_J)(s-s')k'(s-s')\notag\\
    & + |c_I-c_J|k(s-s') \bigr]\dd s'\dd s. \label{61}
\end{align}
\end{lemma}

\begin{proof}
Integration by parts gives
\begin{align*}
    \int_I(k(s-s')-k_I(s'))\dd s &=-\int_I\bigl(s-\frac{a+b}{2}\bigr)k'(s-s')\dd s,\\
    \int_J(k(s-s')-k_J(s))\dd s' &=\int_J\bigl(s'-\frac{a'+b'}{2}\bigr)k'(s-s')\dd s'.
\end{align*}
Put $d_0=\frac{a+b-a'-b'}{2}$. Adding and integrating in both variables yields
\begin{equation}\label{62}
    \iint_{I\times J}(2k-k_I-k_J) =\iint_{I\times J}\bigl[-(s-s')k'(s-s')+d_0k'(s-s')\bigr].
\end{equation}
We show that the integral of the last term is nonpositive. The change of variables $\xi=s-s'$ gives the density
\[
    |I\cap(J+\xi)| =\max\left\{0,\min\left\{|I|,|J|,  \frac{|I|+|J|}{2}-|\xi-d_0|\right\}\right\}.
\]
This expression is nonincreasing in $|\xi-d_0|$. Since
$|\xi+d_0|^2-|\xi-d_0|^2=4\xi d_0$, we have
\[
    d_0\bigl(|I\cap(J+\xi)|-|I\cap(J-\xi)| \bigr) \ge 0 \qquad (\xi>0).
\]
Since $k'$ is odd and nonpositive on $(0,\infty)$,
\begin{align*}
    d_0\iint_{I\times J}k'(s-s')\dd s'\dd s &=\int_0^\infty k'(\xi)d_0\bigl(|I\cap(J+\xi)|-|I\cap(J-\xi)|\bigr)\dd \xi\le0.
\end{align*}
Thus \eqref{62} implies
\[
    \iint_{I\times J}(2k-k_I-k_J)\le-\iint_{I\times J}(s-s')k'(s-s').
\]
By symmetry, assume $c_I\le c_J$. Since $k_J\ge0$,
\begin{align*}
    &\iint_{I\times J}\bigl[c_I(k-k_I)+c_J(k-k_J)\bigr]\\ \quad = & c_I\iint_{I\times J}(2k-k_I-k_J) +(c_J-c_I)\iint_{I\times J}(k-k_J)\\
    \quad \le & -c_I\iint_{I\times J}(s-s')k'(s-s') + (c_J-c_I)\iint_{I\times J}k.
\end{align*}
This proves \eqref{61}.
\end{proof}

\subsection{Directional pressure identities.} We use $\dd S$ to denote surface measure on $\Sph^{n-1}$ and set
\[
    \Ee a(e)=\frac{1}{|\Sph^{n-1}|}\int_{\Sph^{n-1}}a(e)\dd S(e), \qquad \Pe(A)=\frac{1}{|\Sph^{n-1}|}\int_A\dd S(e).
\]
Fix $0\le\psi\in C^\infty([0,\infty))$ such that $\psi>0$ near zero and $\supp\psi\subset[0,\frac{1}{2}]$. For $0<\theta<\frac{1}{4}$ and $r_0<|x|<r_1$, define
\begin{equation}\label{63}
    h_e(x)=\frac{w_0(x)}{Z_\theta} \psi\!\left(\frac{|e-e_r(x)|^2}{\theta^2}\right),\qquad Z_\theta=\Ee\psi\!\left(\frac{|e-e_r(x)|^2}{\theta^2}\right),
\end{equation}
and extend $h_e$ by zero outside the annulus. By rotation invariance, $Z_\theta$ is independent of $x$. Spherical coordinates give
\[
    \begin{aligned} 
        Z_\theta &=\frac{|\Sph^{n-2}|}{|\Sph^{n-1}|}\theta^{n-1}  \int_0^1\psi(s^2)s^{n-2}  \left(1-\frac{\theta^2s^2}{4}\right)^{\frac{n-3}{2}}\dd s,
    \end{aligned}
\]
which yields $c\theta^{n-1}\le Z_\theta\le C\theta^{n-1}$.

Since $w_0$ is compactly supported in the annulus, the zero extension satisfies $h_e\in C_c^2(\R^n)$. The chain rule gives
\begin{equation}\label{64}
 \begin{gathered}
    \Ee h_e=w_0,\qquad \sup_e\|h_e\|_\infty\le C\theta^{1-n}\|w_0\|_\infty,\\
    \sup_e\|\nabla h_e\|_\infty \le C\theta^{-n}(1+r_0^{-1})\|w_0\|_{C^2} \le C Q_0^2\theta^{-n},\\
     \sup_e\|h_e\|_{C^2} \le C\theta^{-(n+1)}(1+r_0^{-1})^2\|w_0\|_{C^2} \le C Q_0^3\theta^{-(n+1)}.
\end{gathered}
\end{equation}
For each $x$ in the annulus, invariance under rotations fixing $e_r$ gives
\[
    \Ee\bigl[h_e(v_e^2-v_r^2)\bigr] =\Ee\bigl[h_e(1-(e\cdot e_r)^2)\bigr] \left(\frac{|\nabla v|^2-v_r^2}{n-1}-v_r^2\right).
\]
On the support of $h_e$,
\[
    0\le1-(e\cdot e_r)^2\le|e-e_r|^2\le\theta^2,\qquad \left|\frac{|\nabla v|^2-v_r^2}{n-1}-v_r^2\right| \le|\nabla v|^2.
\]
Consequently,
\begin{equation}\label{65}
    \left|\Ee\int h_ev_e^2-\int w_0v_r^2\right| \le\theta^2\int w_0|\nabla v|^2 \le\theta^2\|w_0\|_\infty E_\rho.
\end{equation}
For $e\in\Sph^{n-1}$ and $z\ne0$, define
\begin{equation}\label{66}
    P_{\lambda,e}(z)=-(z\cdot e)\partial_eK_\lambda(z) =-K_\lambda'(|z|)\frac{(z\cdot e)^2}{|z|}\ge0.
\end{equation}
For the parallel lines of Lemma~\ref{60}, the identity $-(s-s')k'(s-s')=P_{\lambda,e}(x-y)$ follows directly from the chain rule.

Gaussian differentiation gives
\begin{equation}\label{67}
    P_{\lambda,e}(z)=\int_0^\infty e^{-\lambda t} \frac{(z\cdot e)^2}{2t}\Gamma(z,t)\dd t\ge0.
\end{equation}
Since $\int (z\cdot e)^2\Gamma(z,t)\dd z=2t$,
\[
    \|P_{\lambda,e}\|_1=\int_0^\infty e^{-\lambda t}\dd t=\lambda^{-1}.
\]
For its distributional representation, Fourier differentiation gives
\begin{align*}
    \widehat{P_{\lambda,e}}(\xi) &=e\cdot\nabla_\xi\bigl[(\xi\cdot e)\widehat K_\lambda( \xi)\bigr]
    =(2\pi)^{-\frac{n}{2}}\left[\frac{1}{\lambda+|\xi|^2}  -\frac{2(\xi\cdot e)^2}{(\lambda+|\xi|^2)^2}\right],\\
    \widehat{K_\lambda*K_\lambda}(\xi)
    & =\frac{(2\pi)^{-\frac{n}{2}}}{(\lambda+|\xi|^2)^2}.
\end{align*}
These Fourier transform formulas give
\begin{equation}\label{68}
    P_{\lambda,e}=K_\lambda+2\partial_{ee}(K_\lambda*K_\lambda) \qquad\text{in the sense of distributions}.
\end{equation}

For each fixed direction $e$, put
\begin{equation}\label{69}
    P_{h_e}=\iint h_e(x)P_{\lambda,e}(x-y)g(x)g(y)\dd x\dd y.
\end{equation}
The following identity holds:
\begin{equation}\label{70}
    P_{h_e}=\int h_e gW-2\int h_e W_e^2+R_{h_e},
\end{equation}
where
\begin{align}
    R_{h_e}&=\int(\partial_{ee}h_e)W^2 + 2 \int\bigl(2\nabla h_e\cdot\nabla W+(\Delta h_e) W \bigr)  \partial_{ee} R_\lambda W. \label{71}
\end{align}
Next, we prove that
\begin{align}
    |R_{h_e}| & \le C\|h_e\|_{C^2}  \bigl(\|W\|_2\|\nabla W\|_2+\|W\|_2^2\bigr), \label{72}
\end{align}
where the constant $C$ is uniform for $\lambda\ge1$ and $e\in\Sph^{n-1}$. Since $g\in L^2$ and $h_e\in C_c^2$, we have $W\in H^2$ and $h_eW\in H^2$. The product rule and integrations by parts below are therefore valid in the weak sense. By \eqref{68}, the following identity holds in $L^2$:
\begin{align*}
    P_{\lambda,e}*g &=K_\lambda*g + 2 \partial_{ee} \bigl[(K_\lambda*K_\lambda)* g\bigr]\\
    & = W + 2\partial_{ee}R_\lambda W.
\end{align*}
Self-adjointness and commutation with derivatives yield
\begin{align*}
    P_{h_e} &=\int h_e g(P_{\lambda,e}*g) = \int h_e gW+2\int h_e g\,\partial_{ee}R_\lambda W\\
    &=\int h_e gW+2\int h_e g\,R_\lambda W_{ee}\\
    &=\int h_e gW+2\int R_\lambda(h_e g)W_{ee}.
\end{align*}
The product rule gives
\[
    (-\Delta+\lambda)(h_e W) = h_e g-2\nabla h_e\cdot\nabla W-(\Delta h_e)W,
\]
and hence
\[
    R_\lambda(h_e g)=h_e W+R_\lambda\bigl(2\nabla h_e\cdot\nabla W+(\Delta h_e)W\bigr).
\]
Two integrations by parts in the direction $e$ give
\begin{align*}
    2\int h_e WW_{ee} &=-2\int h_e W_e^2-2\int(\partial_e h_e)WW_e\\
    &=-2\int h_e W_e^2-\int(\partial_e h_e)\partial_e(W^2)\\
    &=-2\int h_e W_e^2+\int(\partial_{ee}h_e)W^2.
\end{align*}
Self-adjointness of $ R_\lambda $ proves \eqref{70} with $R_{h_e}$ as in \eqref{71}. Moreover,
\[
    \|\partial_{ee}R_\lambda W\|_2^2 =\int\frac{|\xi\cdot e|^4}{(\lambda+|\xi|^2)^2} |\widehat W(\xi)|^2\dd\xi \le\|W\|_2^2.
\]
Cauchy--Schwarz gives \eqref{72}.

\begin{lemma}\label{73}
    For every $e\in\Sph^{n-1}$,
    \begin{equation}\label{74}
    \begin{split}
        |R_{h_e}|+2\left| \int h_e( W_e^2 - v_e^2) \right| + \int h_e g | \CC|
        & \le A_\theta \bigl( \lambda^{-\frac{1}{2}} \sqrt{E_\rho D_\rho} + \lambda^{-1} D_\rho\bigr)\\
        & \quad + B_{\theta,\lambda} N_\rho.
    \end{split}
    \end{equation}
    The estimate is uniform in $e\in\Sph^{n-1}$, with
    \[
        A_\theta\le C Q_0^3\theta^{-(n+1)},\qquad B_{\theta,\lambda}\le C Q_0^{n+5}\theta^{-(n+1)}(1+\lambda).
    \]
\end{lemma}

\begin{proof}
    By \eqref{56},
    \begin{align*}
        \|W\|_2\|\nabla W\|_2+\|W\|_2^2 &\le C(1+h^{-4})N_\rho+C\lambda^{-1}\sqrt{E_\rho D_\rho}\\
        &\quad+C(\lambda^{-\frac{3}{2}}+\lambda^{-2})D_\rho.
    \end{align*}
    Moreover,
    \[
        \|\nabla\CC\|_2\le C(1+h^{-2})\|v\|_2 + \frac{\sqrt{D_\rho}}{2\sqrt\lambda}.
    \]
    Using $W_e=v_e+\CC_e$ on $\supp h_e$, we obtain
    \begin{align*}
        |R_{h_e}|+2\left|\int h_e(W_e^2-v_e^2)\right| &\le C Q_0^7 \theta^{-(n+1)} N_\rho\\
        & \quad + C Q_0^3\theta^{-(n+1)} \left(\lambda^{-\frac{1}{2}} \sqrt{E_\rho D_\rho} + \lambda^{-1} D_\rho\right).
    \end{align*}
    Finally, \eqref{59} gives
    \[
        \int h_e g|\CC| \le C Q_0^3\theta^{-(n+1)}\lambda^{-1}D_\rho +C Q_0^{n+5}\theta^{-(n+1)}(1+\lambda)N_\rho.
    \]
    Since $n\ge3$, these inequalities prove \eqref{74}.
\end{proof}

\subsection{Level-set decomposition and comparison.} The following proposition is the main result of this section.

\begin{proposition}\label{75}
    Fix $n\ge3$, $p>1$, $\alpha\ge0$, and $R_*>0$. Let $J\subset\R$ be a fixed bounded open time interval. Let $\rho=r_1+6h\le R_*$, where $r_1>0$ and $0<h\le1$. Let $v\ge0$ be a $C^{2,1}$ solution of \eqref{25} in a neighborhood of $\overline{B_\rho}\times J$, and suppose that
    \[
        \sup_{t\in J}\int_{B_\rho}v(x,t)\dd x\le C_0\rho^{n-\kappa}.
    \]
    Let $ w_0 \in C_c^2(\{r_0<|x|<r_1\}) $ be nonnegative and radial, with $ 0 < r_0 < r_1 $. There are a positive integer $m_p=m_p(n,p)$ and a constant $ C = C(n,p,\alpha,R_*) > 0 $ such that, for every $t\in J$ and $ 0 < \varepsilon \le 1$,
    \begin{equation}\label{76}
        \int w_0 \left(v_r^2-\frac{p-1}{2(p+1)}|x|^\alpha v^{p+1} \right) \le \varepsilon(E_\rho+D_\rho)+C Q^{m_p},
    \end{equation}
    where
    \begin{equation}\label{77}
        Q= 2 + C_0\rho^{n-\kappa}+ \varepsilon^{-1}+ h^{-1}+ r_0^{-1} + \rho^{-1} + \|w_0\|_{C^2(\R^n)}.
    \end{equation}
    Both $C$ and $m_p$ are independent of $J$ and $t$. \end{proposition}

To prove this proposition, we decompose the superlevel sets
of the localized source into intervals along fixed directions. We begin by introducing the necessary notation.

Choose
\begin{equation}\label{78}
    0<\ell\le\min\{\frac{h}{4},\frac{r_0}{4},\frac{1}{4},\frac{1}{4\sqrt\lambda}\}.
\end{equation}
For $0<\tau<g(x)$, define
\[
    I(e,x,\tau) = \left\{ y \in x + \R e: \min_{z\in[x,y]} g(z) > \tau \right\}.
\]
Set
\begin{equation}\label{79}
    d_\pm(e,x,\tau)=\inf\{r\ge0:g(x\pm re)\le\tau\},\qquad x_I^\pm=x\pm d_\pm e.
\end{equation}
Continuity and compact support of $g$ give
\[
    0<d_\pm<\infty,\qquad I=\{x+se:-d_-<s<d_+\},\qquad |I|=d_-+d_+,\qquad g(x_I^\pm)=\tau.
\]
Define the interval weight
\begin{equation}\label{80}
    \chi_I=
    \begin{cases}
        \displaystyle\inf_{y\in I}h_e(y),& |I|\le\ell,\\
        0,& |I|>\ell.
    \end{cases}
\end{equation}
The weight is constant on each component, with $0\le\chi_I\le h_e$. Call $I$ retained if $\chi_I>0$. Then
\[
    |I|\le\ell,\qquad h_e\ge\chi_I>0\ \text{on }\overline I,\qquad \overline I\subset\{r_0<|x|<r_1\}.
\]
Thus $\eta=1$ and $|x|^\alpha$ is smooth and positive on
$\overline I$.

For a function $z$ on the supporting annulus, put $z_I=\frac{z(x_I^-)+z(x_I^+)}{2}$ on intervals with positive weight, and define
\begin{align}
    S_z(e)&=\int\int_0^{g(x)}\chi_I\bigl(z(x)-z_I\bigr)\dd\tau\dd x, \label{81}\\
    L_z(e)&=\int z(x)\int_0^{g(x)}(h_e(x)-\chi_I)\dd\tau\dd x. \label{82}
\end{align}
Integrands are zero outside the supporting annulus, and an endpoint term with zero interval weight is defined to be zero. We use $z=v,W,\CC$. In particular $L_v\ge0$, and all evaluations occur inside $B_{r_1+2h}$.

For $x$ in the annulus and $\tau\ge0$, let $a(x,\tau)\ge0$ solve
\begin{equation}\label{83}
    |x|^\alpha a(x,\tau)^p+\lambda a(x,\tau)=\tau.
\end{equation}
The left-hand side is strictly increasing in $a$. Since $\eta=1$ on the annulus,
\[
    0\le\tau\le g(x)=|x|^\alpha v(x)^p+\lambda v(x) \quad\Longrightarrow\quad 0\le a(x,\tau)\le v(x).
\]
For $\tau>0$, implicit differentiation gives
\[
    (p|x|^\alpha a^{p-1}+\lambda)\nabla a =-\alpha |x|^{\alpha-2}x a^p.
\]
Consequently,
\begin{equation}\label{84}
    |\nabla\log a(x,\tau)| \le \frac{\alpha}{p|x|} \le \frac{\alpha}{pr_0}.
\end{equation}
On a retained interval, $v(x_I^\pm)=a(x_I^\pm,\tau)$. For $0<\tau<g(x)$, the segment from $x$ to either endpoint lies in $\overline I\subset\{r_0<|x|<r_1\}$. Thus \eqref{84} implies
\[
    \left|\log\frac{a(x_I^\pm,\tau)}{a(x,\tau)}\right| \le\int_{[x,x_I^\pm]}|\nabla\log a(y,\tau)|\dd s \le\frac{\alpha|x_I^\pm-x|}{pr_0} \le\frac{\alpha\ell}{pr_0}.
\]
Since $ \ell \le \frac{r_0}{4} $ and $a(x,\tau)\le v(x)$,
\begin{equation}\label{85}
\begin{aligned}
    |v_I-a(x,\tau)| & \le( e^{\frac{\alpha\ell}{pr_0}}-1)a(x,\tau)\\
    & \le \frac{\alpha}{p} e^{\frac{\alpha}{4p}} \frac{ \ell}{ r_0}a(x,\tau) \le C\frac{\ell}{r_0}v(x).
\end{aligned}
\end{equation}
For $\alpha=0$, the endpoint values equal $a(x,\tau)$, so the left-hand side is zero. For $x$ in this annulus and $ z \ge 0 $, define
\begin{equation}\label{86}
    F_\lambda(x,z)=\frac{|x|^\alpha}{p+1}z^{p+1} + \frac{\lambda}{2} z^2.
\end{equation}
The substitution $ \tau = |x|^\alpha z^p + \lambda z $ gives
\begin{align}\label{87}
    \int_0^{g(x)}[v(x)-a(x,\tau)]\dd\tau &=\int_0^{v(x)}(v(x)-z)(p|x|^\alpha z^{p-1}+\lambda)\dd z\notag\\
    & = F_\lambda(x,v(x)).
\end{align}
Using $0 \le a(x,\tau) \le v(x)$ and \eqref{85}, we obtain
\[
\begin{aligned}
    h_e(v-a) & = \chi_I(v-v_I) + \chi_I(v_I-a)+(h_e-\chi_I)(v-a)\\
    & \le \chi_I (v-v_I) + C r_0^{-1} \ell h_ev+(h_e-\chi_I) v.
\end{aligned}
\]
Integration over $\tau$ and $x$ therefore gives
\begin{equation}\label{88}
    \int h_eF_\lambda(x,v) \le S_v(e)+ L_v(e) + C r_0^{-1} \ell \int h_e gv.
\end{equation}
Define
\begin{equation}\label{89}
    M_e = \iint|h_e(x) - h_e(y)| K_\lambda(x-y) g(x) g(y) \dd x\dd y.
\end{equation}

\begin{lemma}\label{90}
For each fixed direction $e$,
\begin{equation}\label{91}
    2S_W(e)\le P_{h_e}+M_e+2L_W(e).
\end{equation}
Consequently,
\begin{align}
    2\int h_e F_\lambda(x,v) \le P_{h_e} + M_e + 4 L_v(e) + 2 \int h_e g|\CC| + 2|S_{\CC}(e)| + C r_0^{-1} \ell \int h_e gv. \label{92}
\end{align}
\end{lemma}

\begin{proof}
    Choose an orthonormal basis whose last vector is $e$. Then $x=z+se$, $y=z'+s'e$, with $z,z'\in e^\perp$, and $\dd x=\dd z\dd s$, $\dd y=\dd z'\dd s'$. Let $\mathcal I(z,\tau)$ denote the components of $\{s:g(z+se)>\tau\}$. For $z\ne z'$, take $I=(a,b)\in\mathcal I(z,\tau)$ and $J=(a',b')\in\mathcal I(z',\sigma)$. Their endpoints are $x_I^-=z+ae$, $x_I^+=z+be$, $y_J^-=z'+a'e$, and $y_J^+=z'+b'e$. Lemma~\ref{60}, with $c_I=\chi_I$ and $c_J=\chi_J$, gives
    \begin{align*}
        & \int_a^b \int_{a'}^{b'}(\chi_I+\chi_J) K_\lambda(x-y) \dd s' \dd s\\
        & \quad \le \int_a^b \int_{a'}^{b'} \Bigl[ \frac{\chi_I}{2} \bigl( K_\lambda( x_I^--y) + K_\lambda(x_I^+-y) \bigr)\\
        & \qquad\qquad+\frac{\chi_J}{2} \bigl( K_\lambda(x-y_J^-) + K_\lambda(x  -y_J^+) \bigr)\\
        & \qquad \qquad + \min(\chi_I,\chi_J) P_{\lambda,e}(x-y) + |\chi_I-\chi_J|K_\lambda(x-y)\Bigr]\dd s'\dd s.
    \end{align*}
    For each fixed $z\in e^\perp$ and $\tau>0$,
    \begin{align*}
        \sum_{I \in \mathcal I(z,\tau)} \mathbf1_I(s) &=\mathbf1_{\{\tau<g(z+se)\}},\\
        \sum_{I \in \mathcal I(z,\tau)} \chi_I\mathbf1_I(s) & = \chi_{I(e,z+se,\tau)} \mathbf1_{\{\tau<g(z+se)\}}.
    \end{align*}
    The expression $\chi_{I(e,x,\tau)}$ is extended by zero for $ \tau \ge g(x)$. In the following sums, $ I \in \mathcal I(z,\tau) $ and $ J \in \mathcal I(z',\sigma) $. The first term on the left is
    \begin{align*}
        & \sum_{I,J} \int_I \int_J \chi_I K_\lambda(x-y)\dd s'\dd s\\
        &\quad= \int_\R \int_\R \left( \sum_I \chi_I \mathbf1_I(s) \right) \left( \sum_J \mathbf1_J(s') \right) K_\lambda(x-y) \dd s'\dd s.
    \end{align*}
    For fixed $x$ and $I$,
    \begin{align*}
        & \int_{e^\perp} \int_0^\infty \int_\R \left( \sum_{ J \in \mathcal I(z', \sigma)} \mathbf1_J (s') \right) K_\lambda(x-z'-s'e) \dd s'\dd\sigma \dd z'\\
        & \quad = \int_{\R^n} K_\lambda(x-y) \int_0^\infty \mathbf1_{\{\sigma<g(y)\}} \dd \sigma \dd y\\
        & \quad = \int_{\R^n} K_\lambda(x-y) g(y) \dd y = W(x).
    \end{align*}
    The remaining integrations therefore give
    \begin{align*}
        & \int_{e^\perp} \int_0^\infty \sum_{I \in \mathcal I(z,\tau)} \int_\R \chi_I\mathbf1_I(s)W(z+se)\dd s\dd\tau\dd z\\
        & = \int_{\R^n} W(x) \int_0^\infty \chi_{I(e,x,\tau) }\mathbf1_{\{ \tau < g(x)\}} \dd \tau \dd x\\
        & = \int_{\R^n} W(x) \int_0^{g(x)} \chi_I \dd\tau\dd x.
    \end{align*}
    By symmetry, the term with $\chi_J$ gives the same quantity. For the endpoint term, fix $x_I^\pm$ and use
    \[
        \frac{1}{2}\int_{\R^n}  \bigl(K_\lambda(x_I^--y)+K_\lambda(x_I^+-y)\bigr)g(y)\dd y  =\frac{W(x_I^-)+W(x_I^+)}{2}.
    \]
    Thus the first endpoint term contributes
    \begin{align*}
        &\int_{e^\perp}\int_0^\infty \sum_{I\in\mathcal I(z,\tau)}\int_I \chi_I\frac{W(x_I^-)+W(x_I^+)}{2}\dd s\dd\tau\dd z\\
        &=\int_{\R^n}\int_0^{g(x)}  \chi_I\frac{W(x_I^-)+W(x_I^+)}{2}\dd\tau\dd x,
    \end{align*}
    and the second contributes the same integral. No change of variables involving the endpoints is made. For the pressure term, positivity and $\min(\chi_I,\chi_J)\le h_e(x)$ give
    \begin{align*}
        &\iint P_{\lambda,e}(x-y) \int_0^{g(x)} \int_0^{g(y)} \min(\chi_I,\chi_J) \dd \sigma \dd \tau \dd x \dd y\\
        & \quad \le \iint h_e(x) P_{\lambda,e}(x-y) g(x)g(y)\dd x\dd y =P_{h_e}.
    \end{align*}
    Finally,
    \[
        |\chi_I-\chi_J|\le|h_e(x)-h_e(y)|+( h_e(x)-\chi_I)+(h_e(y)-\chi_J),
    \]
    so symmetry yields
    \begin{align*}
        &\iint K_\lambda(x-y)\int_0^{g(x)}\int_0^{g(y)} | \chi_I - \chi_J| \dd \sigma \dd \tau \dd x \dd y\\
        & \quad \le M_e+2\int\int_0^{g(x)}(h_e(x)-\chi_I) \left( \int K_\lambda(x-y) g(y) \dd y \right) \dd \tau \dd x\\
        & \quad = M_e + 2 \int W(x) \int_0^{g(x)} (h_e(x)-\chi_I) \dd \tau \dd x = M_e + 2 L_W(e).
    \end{align*}
    Consequently,
    \[
    \begin{aligned}
        2 \int_{\R^n}\int_0^{g(x)} \chi_I W(x) \dd \tau \dd x & \le 2 \int_{\R^n} \int_0^{g(x)} \chi_I W_I \dd \tau \dd x\\
        & \quad + P_{h_e} + M_e + 2L_W(e).
    \end{aligned}
    \]
    Subtracting the first term on the right proves \eqref{91}.

    On retained intervals, $W=v+\CC$. The linearity of endpoint averages gives $S_W=S_v+S_{\CC}$. Also
    \[
        L_W-L_v=\int\CC(x)\int_0^{g(x)}(h_e(x)-\chi_I) \dd\tau\dd x \le\int h_eg|\CC|.
    \]
    Multiplying \eqref{88} by $2$ and using these relations gives
    \begin{align*}
        2\int h_eF_\lambda &\le2S_v+2L_v+Cr_0^{-1}\ell\int h_egv\\
        & =2 S_W - 2 S_{\CC}+2L_v + Cr_0^{-1}\ell\int h_e g v \\
        & \le P_{h_e} + M_e + 2L_W + 2L_v + 2|S_{\CC}| + Cr_0^{-1} \ell \int h_e g v\\
        & \le P_{h_e} + M_e + 4 L_v + 2 \int h_e g|\CC| + 2 |S_{\CC}| + C r_0^{-1} \ell \int h_e gv.
    \end{align*}
    This proves \eqref{92}.
\end{proof}

\begin{remark}\label{93}
    Since $g\in C_c(\R^n)$ and $K_\lambda\in L^1$, $W$ is continuous and
    \[
        \|W\|_\infty\le\lambda^{-1}\|g\|_\infty.
    \]
    Using $0\le\chi_I\le h_e\le C\theta^{1-n}\|w_0\|_\infty$ and $\|P_{\lambda,e}\|_1 = \lambda^{-1} $, we obtain
    \[
        \begin{aligned}
            & \int_{\R^n} \int_0^{g(x)} \chi_I \bigl(W(x) + W_I\bigr) \dd \tau \dd x + L_W(e) + P_{h_e} + M_e\\
            & \qquad \le C \theta^{1-n} \| w_0 \|_\infty \lambda^{-1} \|g\|_\infty \|g\|_1<\infty.
        \end{aligned}
    \]
    The other layer integrals are finite because $v,\CC$ are locally bounded. These bounds justify the changes in the order of integration and the subtractions above. They are used only to establish finiteness: $\|g\|_\infty$ does not enter the quantitative constants.
\end{remark}

\subsection{Error estimates and proof of Proposition \ref{75}.} We estimate $L_v(e)$, $S_{\CC}(e)$, and $M_e$ in turn. The other correction terms are controlled by Lemma~\ref{73}.

First, we estimate $L_v(e)$.

\begin{lemma}\label{94}
    There is $C_*>0$, depending only on $n,p,\alpha,R_*$, such that, if
    \begin{equation}\label{95}
        \tau_*\ge C_*\lambda^{\frac{p}{p-1}},
    \end{equation}
    then
    \begin{equation}\label{96}
    \begin{aligned}
        \Ee L_v(e) & \le C Q_0^2 \theta^{-n} ( \ell + \ell^{-n} \tau_*^{-\frac{1}{p}} ) \int_{ B_{r_1+2h} } g v\\
        & \quad + C \theta^{1-n} \| w_0 \|_\infty \tau_*C_0 \rho^{n-\kappa}.
    \end{aligned}
    \end{equation}
\end{lemma}

\begin{proof}
    For $|I|\le\ell$, including the case $\chi_I=0$,
    \[
        0 \le h_e(x)- \chi_I \le \|\nabla h_e \|_\infty\ell \le C Q_0^2 \theta^{-n}\ell.
    \]
    Thus
    \[
        \Ee \int_{B_{r_1}} v(x) \int_0^{g(x)} (h_e(x) - \chi_I) \mathbf1_{\{|I|\le\ell\}} \dd \tau \dd x \le C Q_0^2 \theta^{-n} \ell \int_{B_{r_1}} g v.
    \]
    If $|I| > \ell$ and $ x \in B_{r_1} $, one of $ d_+, d_- $ exceeds $ \frac{\ell}{2} $. In that direction the segment $0<r<\frac{\ell}{2}$ lies in $\{g>\tau\}$ and is contained in $B_{r_1+2h}$ by $\ell\le \frac{h}{4}$ in \eqref{78}. There $g=|x|^\alpha v^p+\lambda v$ and $|x|^\alpha\le R_*^\alpha$. If $g(y)>\tau$, at least one of $R_*^\alpha v(y)^p>\frac{\tau}{2}$ and $\lambda v(y)>\frac{\tau}{2}$ holds. Therefore
    \begin{equation}\label{97}
        v(y) \ge \min \left\{\left(\frac{\tau}{2R_*^\alpha} \right)^{\frac{1}{p}}, \frac{\tau}{2\lambda} \right\}.
    \end{equation}
    Taking $ C_*: = \max\{1, 2R_*^{-\frac{\alpha}{p-1}}\}$, we obtain $v(y)\ge c\tau^{\frac{1}{p}}$ whenever $ g(y) > \tau $ and $ \tau \ge C_* \lambda^{\frac{p}{p-1}} $. Fix $x\in B_{r_1}$ and $ C_*\lambda^{\frac{p}{p-1}} \le \tau < g(x) $, and put
    \[
        E_\pm = \{e \in \Sph^{n-1 }: d_\pm(e,x,\tau) > \frac{\ell}{2}\}.
    \]
    Since $B_{\frac{\ell}{2}}(x)\subset B_{r_1+2h}$, polar coordinates give
    \begin{align*}
        C_0 \rho^{n-\kappa} & \ge \int_{B_{\frac{\ell}{2}}(x)} v(y) \dd y\\
        & \ge c \tau^{\frac{1}{p}} \int_{E_+} \int_0^{\frac{\ell}{2}} r^{n-1} \dd r \dd S(e) = \frac{c}{n} \tau^{\frac{1}{p}} (\frac{\ell}{2})^n |E_+|.
    \end{align*}
    Moreover,
    \[
        E_-=-E_+,\qquad \{|I(e,x,\tau)|>\ell\}\subset E_+\cup E_-.
    \]
    Hence
    \begin{equation}\label{98}
    \begin{aligned}
        \Pe\{|I(e,x,\tau)| > \ell\} &\le\frac{|E_+|+|E_-|}{|\Sph^{n-1}|} = \frac{2|E_+|}{|\Sph^{n-1}|}\\
         &\le C C_0\rho^{n-\kappa}\ell^{-n}\tau^{-\frac{1}{p}}.
    \end{aligned}
    \end{equation}
    For the long intervals with $\tau<\tau_*$,
    \[
        \Ee \int_{B_{r_1}} v(x) \int_0^{\min\{\tau_*,g(x)\}} h_e(x) \mathbf1_{\{|I|>\ell\}}\dd\tau\dd x \le C\theta^{1-n}\|w_0\|_\infty \tau_* C_0\rho^{n-\kappa}.
    \]
    For $\tau\ge\tau_*$, \eqref{98} gives
    \begin{align*}
        & \Ee \int_{B_{r_1}\cap\{g > \tau_*\}} v(x) \int_{\tau_*}^{g(x)} h_e(x) \mathbf1_{\{|I|>\ell\}}\dd\tau\dd x\\
        & \quad \le C \theta^{1-n} \|w_0\|_\infty C_0 \rho^{n-\kappa} \ell^{-n}
       \int_{B_{r_1} \cap\{ g > \tau_*\}} v(x)  \int_{\tau_*}^{g(x)} \tau^{-\frac{1}{p}}\dd\tau\dd x \\
        & \quad \le C Q_0^2 \theta^{1-n}  \ell^{-n} \tau_*^{-\frac{1}{p}} \int_{B_{r_1}} g v.
    \end{align*}
    Adding these three estimates proves \eqref{96}.
\end{proof}

Next, we estimate $S_{\CC}(e)$. Since $p>1$ and $r_0>0$, $g$ is locally Lipschitz on $\{r_0<|x|<r_1\}$. Interior elliptic regularity gives
\[
     W \in C_{\rm loc}^{2,\beta} \quad( 0 < \beta < 1),\qquad \CC = W - v \in C^2,\qquad H=\CC-B\in C^2
\]
on this annulus. Thus the interval integrations by parts are valid. The quantitative estimates use only \eqref{55} and \eqref{54}.

For $q\in C^2([a,b])$, two integrations by parts with the kernel $\frac{1}{2}(s-a)(b-s)$ give
\begin{equation}\label{99}
    \int_a^b\left(q(s)-\frac{q(a)+q(b)}{2}\right)\dd s =-\frac{1}{2}\int_a^b(s-a)(b-s)q''(s)\dd s.
\end{equation}

Since $\chi_I$ is constant on each interval, no derivative of $\chi_I$ appears. Define
\begin{equation}\label{100}
    T_e(x)=\frac{1}{2}\int_0^{g(x)}\chi_Id_+(e,x,\tau)d_-(e,x,\tau)\dd\tau, \qquad \overline T(x)=\Ee T_e(x).
\end{equation}
These densities are nonnegative and supported in the annulus. Since $d_+d_-\le\frac{\ell^2}{4}$ when $\chi_I>0$,
\[
    0\le T_e(x)\le\frac{\ell^2}{8}h_e(x)g(x),\qquad \int T_e|\CC_{ee}| \le\frac{\ell^2}{8}\int h_eg|\CC_{ee}|<\infty.
\]
This bound is used only to justify changes in the order of integration. For $x=z+se$ in $I=(a,b)$, $d_-=s-a$ and $d_+=b-s$. Applying \eqref{99} on each interval gives
\begin{align*}
    S_{\CC}(e) & = \int_{e^\perp}\int_0^\infty \sum_{ \substack{I=(a,b) \in \mathcal I(z,\tau) \\ \chi_I>0}}  \chi_I \int_a^b\bigl(\CC(z+se)-\CC_I\bigr)\dd s\dd\tau\dd z\\
    & = - \frac{1}{2} \int_{e^\perp} \int_0^\infty \sum_{ \substack{ I = (a,b) \in \mathcal I(z,\tau) \\ \chi_I>0}} \chi_I \int_a^b(s-a)(b-s) \CC_{ee}(z+se)\dd s \dd\tau\dd z\\
    & = -\int_{\R^n} \left(\frac{1}{2} \int_0^{g(x)} \chi_I d_+(e,x,\tau) d_-(e,x,\tau)\dd\tau\right) \CC_{ee}(x) \dd x.
\end{align*}
Thus
\begin{equation}\label{101}
    S_{\CC}(e)=-\int T_e\CC_{ee}.
\end{equation}

\begin{lemma}\label{102}
    Under \eqref{78},
    \begin{align}
        \overline T(x) & \le C \theta^{1-n} \|w_0\|_\infty  \int_{|y-x|<\ell} |x-y|^{2-n}g(y)\dd y \le C\theta^{1-n}\|w_0\|_\infty W(x),\label{103}\\
        \|\overline T\|_2 &\le C\theta^{1-n}\| w_0 \|_\infty [( 1+h^{-2}) \|v\|_2 +\lambda^{-1} \sqrt{D_\rho}] ,\label{104}\\
        \| \overline T \|_1 & \le C \theta^{1-n} \| w_0 \|_\infty [( 1+ \lambda^{-1} h^{-2}) C_0 \rho^{n-\kappa} + \lambda^{-1} \sqrt{D_\rho}]. \label{105}
\end{align}
The constants are independent of $\lambda,\ell$ and the solution.
\end{lemma}

\begin{proof}
    Fix $e,x$ and write $I=I(e,x,\tau)$. For $0\le r<\ell$,
    \begin{align*}
        \int_0^{g(x)} \mathbf1_{\{|I|\le\ell\}} \mathbf1_{\{d_\pm>r\}} \dd\tau
        & \le\int_0^\infty \mathbf1_{\{\tau<g(x\pm re)\}} \dd\tau = g (x\pm re).
    \end{align*}
    Using $d^2=2\int_0^\ell r\mathbf1_{\{r<d\}}\dd r$ for $0\le d\le\ell$, Tonelli's theorem gives
    \begin{align}
        \int_0^{g(x)}d_\pm^2\mathbf1_{\{|I|\le\ell\}}\dd\tau & = 2 \int_0^\ell r\int_0^{g(x) } \mathbf1_{\{|I|\le\ell\}} \mathbf1_{\{d_\pm>r\}}\dd\tau\dd r\notag\\ & \le 2 \int_0^\ell r g(x\pm re)\dd r. \label{106}
    \end{align}
    Since $ 2 d_+ d_- \le d_+^2+d_-^2$ and $\chi_I\le C\theta^{1-n} \|w_0\|_\infty\mathbf1_{\{|I|\le\ell\}}$, we have
    \begin{align*}
        \overline T(x) & = \frac{1}{2} \Ee \int_0^{g(x)} \chi_Id_ + d_- \dd\tau  \le \frac{C\theta^{1-n} \|w_0\|_\infty}{4} \Ee \int_0^{g(x)} (d_+^2 + d_-^2) \mathbf1_{\{|I|\le\ell\}} \dd \tau\\
        & \le \frac{C\theta^{1-n}\|w_0\|_\infty}{2}\Ee\int_0^\ell r \bigl( g(x+re) + g(x-re) \bigr) \dd r \\
        & = \frac{C \theta^{1-n} \|w_0\|_\infty}{ |\Sph^{n-1}| } \int_{\Sph^{n-1}} \int_0^\ell r g(x+re) \dd r \dd S(e)\\
        & = \frac{C\theta^{1-n} \|w_0\|_\infty}{ |\Sph^{n-1}| } \int_{|y-x|<\ell} |x-y|^{2-n} g(y) \dd y\\
        & \le C \theta^{1-n} \|w_0\|_\infty \int_{|y-x|<\ell} K_\lambda(x-y) g(y) \dd y \le C\theta^{1-n}\|w_0\|_\infty W(x).
    \end{align*}
    The last line uses \eqref{50} and $\ell\le\frac{1}{4\sqrt\lambda}$. The $L^2$ bound follows from \eqref{56}. For $L^1$, use
    \[
        \|\overline T\|_1 \le C\theta^{1-n}\|w_0\|_\infty\int W =C\theta^{1-n}\|w_0\|_\infty\lambda^{-1}\int g.
    \]
    Equation~\eqref{57} and $\rho\le R_*$ prove \eqref{105}.
\end{proof}

\begin{lemma}\label{107}
    Under \eqref{78},
    \begin{equation}\label{108}
    \begin{aligned}
        \Ee|S_{\CC}(e)| & \le C \theta^{1-n} \|w_0\|_\infty  \Bigl[ ( 1+h^{-2}) \|v\|_2 \sqrt{D_\rho} + \lambda^{-1}D_\rho\\
         & \qquad + h^{-n-2} C_0 \rho^{n-\kappa} \bigl(( 1 + \lambda^{-1} h^{-2}) C_0 \rho^{n-\kappa} + \lambda^{-1} \sqrt{D_\rho} \bigr) \Bigr].
    \end{aligned}
    \end{equation}
\end{lemma}

\begin{proof}
    Use \eqref{101}, $\CC=H+B$, and $T_e\ge0$. For every $e$, $|H_{ee}|\le|D^2H|$ and $|B_{ee}|\le\|D^2B\|_{L^\infty(B_{r_1+2h})}$ on the support. Thus
    \begin{align*}
        \Ee|S_{\CC}(e)| & \le \Ee\int T_e|H_{ee}|+\Ee\int T_e|B_{ee}| \\
        & \le \int \overline  T |D^2H| + \|D^2B\|_{L^\infty(B_{r_1+2h})} \int\overline T\\
        & \le \|\overline T \|_2 \|D^2H\|_2 + \|D^2B \|_{L^\infty(B_{r_1+2h})} \|\overline T\|_1.
\end{align*}
Now apply \eqref{55}, \eqref{54}, \eqref{104}, and \eqref{105}.
\end{proof}

The absolute value is taken for each direction before averaging. Although $\|D^2H\|_2\le\sqrt{D_\rho}$ has no small factor, the comparison of the averaged density with $W$ supplies $\lambda^{-1}$ in front of $D_\rho$.

Finally, we estimate $M_e$. Write $ M_e = M_e^{(1)} + M_e^{(2)} $, where
\begin{align*}
    M_e^{(1)} & = \iint_{|x-y|<\ell} |h_e(x)-h_e(y)| K_\lambda(x-y) g(x) g(y) \dd x\dd y,\\
    M_e^{(2)} & = \iint_{|x-y|\ge\ell} |h_e(x)-h_e(y)| K_\lambda(x-y) g(x) g(y) \dd x\dd y.
\end{align*}
For $M_e^{(1)}$, a nonzero weight difference requires at least one point in $B_{r_1}$. If $|x-y|<\ell$, the other lies in $B_{r_1+2h}$ because $ \ell \le \frac{h}{4}$ by \eqref{78}. The Lipschitz bound for $h_e$ gives
\[
    \mathbf1_{\{|x-y|<\ell\}}|h_e(x)-h_e(y)| \le C Q_0^2 \theta^{-n} \ell \,\mathbf1_{B_{r_1+2h}}(x) \mathbf1_{B_{r_1+2h}}(y).
\]
Hence
\begin{equation}\label{109}
\begin{aligned}
    M_e^{(1)} & \le C Q_0^2 \theta^{-n} \ell \int_{B_{r_1+2h}} g(x) \int_{B_{r_1+2h}} K_\lambda(x-y) g(y) \dd y \dd x\\
    & \le C Q_0^2 \theta^{-n} \ell \int_{B_{r_1+2h}} g W \\
    & \le C Q_0^2 \theta^{-n} \ell \left( E_\rho + \lambda \|v\|_2^2 + \int_{B_{r_1+2h}} g|\CC|\right).
\end{aligned}
\end{equation}
The last term is controlled by \eqref{59}.

For $M_e^{(2)}$, choose $\zeta\in C_c^\infty(B_{r_1+2h})$, $0\le\zeta\le1$, with $\zeta=1$ on $B_{r_1}$ and $\|D^j\zeta\|_\infty\le C_jh^{-j}$ for $1\le j\le4$, and fix a smooth cutoff $0\le\beta\le1$ on $[0,\infty)$, equal to zero on $[0,\frac{1}{2}]$ and one on $[1,\infty)$. With a sufficiently large fixed $C$, define
\begin{equation}\label{110}
    \KK(x,y) = C\theta^{1-n} \|w_0\|_\infty( \zeta(x)+\zeta(y)) K_\lambda(x-y) \beta(\frac{|x-y|}{\ell}).
\end{equation}
By the choices of $\zeta$, $\beta$, and $C$,
\[
    \mathbf1_{\{|x-y|\ge\ell\}} |h_e(x)-h_e(y)|K_\lambda(x-y) \le \KK(x,y).
\]
Therefore
\begin{equation}\label{111}
    M_e^{(2)}\le\iint\KK(x,y)g(x)g(y)\dd x\dd y.
\end{equation}
This kernel is symmetric, nonnegative, and is zero near the diagonal. For fixed $\theta,\lambda,\ell$, its derivatives through order four are bounded on $B_{r_1+5h} \times B_{r_1+5h}$. The same symbol $\KK$ will denote the integral operator with this kernel. Its operator norm on $L^2(\R^n)$ satisfies
\begin{equation}\label{112}
    \|\KK\|_{L^2\to L^2} \le C\theta^{1-n} \|w_0\|_\infty \lambda^{-1}.
\end{equation}
The constant in this estimate is independent of $\ell$. Indeed, the kernel mass gives
\[
    \sup_x \int \KK(x,y) \dd y\le C \theta^{1-n} \|w_0\|_\infty\lambda^{-1}.
\]
Symmetry and the Cauchy--Schwarz inequality give
\begin{align*}
    \|\KK q\|_2^2 & \le \int\left( \int\KK(x,y) \dd y \right) \left( \int\KK(x,y)|q(y)|^2 \dd y \right)\dd x \\
    & \le C \theta^{1-n} \|w_0\|_\infty \lambda^{-1}\int |q(y)|^2 \left( \int\KK(x,y) \dd x \right) \dd y\\
    & \le C \theta^{2(1-n)} \|w_0\|_\infty^2 \lambda^{-2} \|q\|_2^2.
\end{align*}
This proves \eqref{112} with a constant independent of $\ell$.

\begin{lemma}\label{113}
    For every $e\in\Sph^{n-1}$,
    \begin{equation}\label{114}
    \begin{aligned}
        M_e^{(2)} & \le C\theta^{1-n} \|w_0\|_\infty \lambda^{-1} D_\rho\\
        & \quad + C Q_0^7 \theta^{1-n} (1+\lambda)^2 \ell^{-(n+2)}(1+\sqrt{D_\rho}).
    \end{aligned}
    \end{equation}
\end{lemma}

\begin{proof}
    Use the multi-index notation $|\gamma|=\sum_j\gamma_j$ and $\partial^\gamma = \partial_1^{\gamma_1} \cdots \partial_n^{\gamma_n}$. Equation \eqref{52} gives, in distributions,
    \begin{equation}\label{115}
        g = \eta v_t + \sum_{|\gamma|\le2} \partial^\gamma F_\gamma,
    \end{equation}
    where the nonzero coefficients are
    \[
        F_0=\lambda\eta v-v\Delta\eta,\qquad F_{e_j} = 2v \partial_j \eta,\qquad F_{2e_j} = -\eta v \quad( 1\le j\le n).
    \]
    Here $e_j$ is the $j$th unit multi-index; all other coefficients are zero. Each $F_\gamma$ is compactly supported in $B_{r_1+5h}$, and
    \begin{equation}\label{116}
        \sum_{|\gamma|\le2}\|F_\gamma\|_1 \le C(\lambda+h^{-2})C_0\rho^{n-\kappa}.
    \end{equation}
    For $\BB(a,b)=\iint\KK(x,y)a(x)b(y)\dd x\dd y$, symmetry gives
    \begin{equation}\label{117}
    \begin{aligned}
        \BB(g,g) &=\BB(g-\eta v_t,g-\eta v_t)\\
        & \quad + 2\BB(g-\eta v_t,\eta v_t)+\BB(\eta v_t, \eta v_t).
    \end{aligned}
    \end{equation}
    The last term satisfies
    \begin{equation}\label{118}
        |\BB(\eta v_t,\eta v_t)| \le C\theta^{1-n} \|w_0 \|_\infty \lambda^{-1} D_\rho.
    \end{equation}
    Put
    \[
        J_j = \max_{|\gamma|+|\mu|\le j}
       \| \partial_x^\gamma \partial_y^\mu \KK\|_{L^\infty (B_{r_1+5h} \times B_{r_1+5h})}, \qquad j=2,4.
    \]
    On the support of the differentiated kernel, $|x-y| \ge \frac{\ell}{2}$. If $k$ derivatives fall on $\zeta$, their contribution is bounded by
    \[
        h^{-k}\ell^{2-n-j+k} = \left( \frac{\ell}{h} \right)^k \ell^{2-n-j} \le \ell^{2-n-j},
    \]
    because $\ell\le\frac{h}{4}$. Thus \eqref{46} and the cutoff bounds give
    \[
        J_j \le C \theta^{1-n} \|w_0\|_\infty \ell^{2-n-j},\qquad j=2,4.
    \]
    In particular,
    \[
        J_2 \le C Q_0\theta^{1-n}\ell^{-n},\qquad J_4 \le C Q_0 \theta^{1-n} \ell^{-(n+2)}.
    \]
    Integration by parts gives
    \begin{align*}
        \BB(g-\eta v_t,\eta v_t) & = \sum_{|\gamma|\le2}(-1)^{|\gamma|} \iint F_\gamma(x) (\eta v_t)(y)\partial_x^\gamma\KK(x,y)\dd x\dd y,\\
        \BB(g-\eta v_t,g-\eta v_t) & = \sum_{|\gamma|, |\mu|\le2}(-1)^{ |\gamma|+|\mu|} \iint F_\gamma(x)F_\mu(y) \partial_x^\gamma \partial_y^\mu \KK(x,y)\dd x\dd y.
    \end{align*}
    Since $\|\eta v_t\|_1\le C(R_*)\sqrt{D_\rho}$,
    \begin{align*}
        |\BB(g-\eta v_t,\eta v_t)| & \le C J_2(\lambda+h^{-2}) C_0\rho^{n-\kappa} \sqrt{D_\rho},\\
        |\BB(g-\eta v_t,g-\eta v_t)| & \le C J_4( \lambda+h^{-2})^2 C_0^2\rho^{2(n-\kappa)}.
    \end{align*}
    There are no boundary terms, since the coefficients have compact support. Using \eqref{111} and \eqref{117}, we obtain
    \begin{align*}
        M_e^{(2)} & \le |\BB(g-\eta v_t,g-\eta v_t)| + 2|\BB(g-\eta v_t, \eta v_t)|\\
        & \quad + |\BB(\eta v_t,\eta v_t)|\\
         & \le C Q_0^7 \theta^{1-n} (1+\lambda)^2\ell^{-(n+2)} (1+\sqrt{D_\rho}) + C\theta^{1-n}\|w_0\|_\infty\lambda^{-1}D_\rho.
    \end{align*}
    Here $J_4$ contributes $Q_0$, and each factor $(\lambda+h^{-2}) C_0\rho^{n-\kappa}$ contributes at most $C Q_0^3(1+\lambda)$. Their product is bounded by $C Q_0^7(1+\lambda)^2$.
\end{proof}

The coefficient of $D_\rho$ is independent of $\ell$. We are ready to prove Proposition \ref{75}.

\begin{proof}[Proof of Proposition~\ref{75}]
On the support of $h_e$, $\eta=1$ and
\begin{equation}\label{119}
    gv-2F_\lambda(x,v)=\frac{p-1}{p+1}|x|^\alpha v^{p+1}.
\end{equation}
Substitute \eqref{70} into \eqref{92} and use $W=v+\CC$. This gives
\begin{equation}\label{120}
\begin{aligned}
    2 \int h_e \left( v_e^2 - \frac{p-1}{2(p+1)} |x|^\alpha v^{p+1} \right)  & \le |R_{h_e}| + 2\left| \int h_e( W_e^2-v_e^2 ) \right|\\
    & \quad + 3 \int h_e g |\CC| + M_e + 4L_v(e)\\
    & \quad + 2|S_{\CC}(e)| + Cr_0^{-1} \ell\int h_e gv.
\end{aligned}
\end{equation}
Divide by two and average in $e$. Apply \eqref{65}, Lemma~\ref{73}, \eqref{59}, \eqref{96}, \eqref{108}, \eqref{109}, and \eqref{114}. Since
\[
    \int_{B_{r_1+2h}}gv\le E_\rho+\lambda\|v\|_2^2,
\]
the lower-order terms are bounded by multiples of $N_\rho$. The coefficients of $\lambda^{-1}D_\rho$ in the endpoint and far-field estimates are independent of $\ell,\tau_*$. Collecting the estimates gives
\begin{equation}\label{121}
\begin{aligned}
    &\int w_0\left(v_r^2-\frac{p-1}{2(p+1)}|x|^\alpha v^{p+1}\right)\\
    &\quad\le[A_0\theta^2+A_\theta(\ell+\ell^{-n}\tau_*^{-\frac{1}{p}})]E_\rho\\
    &\qquad\quad+A_\theta(\lambda^{-\frac{1}{2}}\sqrt{E_\rho D_\rho} + \lambda^{-1} D_\rho) + B_{\theta,\lambda,\ell,\tau_*} N_\rho.
\end{aligned}
\end{equation}
For a sufficiently large fixed $C=C(n,p,\alpha,R_*)$, we may take
\[
    A_0=CQ_0,\qquad A_\theta=CQ_0^3\theta^{-(n+1)},
\]
and
\[
    B_{\theta,\lambda,\ell,\tau_*} \le C Q_0^{n+7}\theta^{-(n+1)} \bigl[(1+\lambda)^2\ell^{-(n+2)}+1+\tau_*\bigr].
\]
Only the low source layers contribute a positive power of $\tau_*$. The coefficient $A_\theta$ is independent of $\lambda,\ell,\tau_*$. For fixed $0<\varepsilon\le1$, first choose $\theta$, which determines $A_\theta$, and then choose $\lambda$, $\ell$, and $\tau_*$ successively:
\begin{equation}\label{122}
    \left\{ \begin{aligned}
        \theta & =\min\left\{\frac{1}{8}, \sqrt{ \frac{\varepsilon}{64(1+A_0)}} \right\},\\
        \lambda & = \left(1+\frac{64(1+A_\theta)}{\varepsilon} \right)^2,\\
        \ell^{-1} & = 8 + 8 h^{-1} + 8r_0^{-1} + 8 \sqrt \lambda + \frac{64(1+A_\theta)}{\varepsilon},\\
        \tau_* & = C_* \lambda^{\frac{p}{p-1}} + \left(\frac{64(1+A_\theta)\ell^{-n}}{\varepsilon}\right)^p.
    \end{aligned} \right.
\end{equation}
Each choice uses only the fixed data and the preceding parameters. These choices satisfy \eqref{78} and \eqref{95}. Using $\sqrt{E_\rho D_\rho}\le\frac{E_\rho+D_\rho}{2}$, the energy and dissipation terms in \eqref{121} are bounded by $\frac{\varepsilon}{2}(E_\rho+D_\rho)$.

Since $ Q  = Q_0 + \varepsilon^{-1} $, the parameter choices give
\[
    \begin{aligned}
        \theta^{-1} & \le CQ,\qquad \sqrt\lambda+\ell^{-1}\le CQ^{n+5},\\
        \tau_* & \le C Q^{(n+5) \max\{p(n+1), \frac{2p}{p-1}\}}.
    \end{aligned}
\]
Consequently,
\[
\begin{aligned}
    B_{\theta,\lambda,\ell,\tau_*} & \le C Q^{2n+8} \left[ Q^{(n+5)(n+6)} + Q^{(n+5) \max\{p(n+1),\frac{2p}{p-1}\}}\right]\\
    & \le C Q^{2n+8+(n+5) \max\{n+6,p(n+1), \frac{2p}{p-1} \}}.
\end{aligned}
\]
Set
\[
    m_p=\left\lceil (n+2)\left[2n+8+(n+5) \max \left\{n+6,p(n+1), \frac{2p}{p-1}\right\}\right] + n + 3\right\rceil.
\]
Apply Young's inequality term by term as in the proof of Lemma~\ref{27} with $ L = B_{\theta,\lambda,\ell,\tau_*} $ and absorption parameter $ \frac{\varepsilon}{2} $. The remainder with the largest power of $ Q $ is
\[
    C \varepsilon^{-(n+1)} B_{\theta,\lambda,\ell, \tau_*}^{n+2} C_0^2\rho^{2(n-\kappa)}\le C Q^{m_p}.
\]
The other remainders have no larger power of $Q$, and hence
\[
    B_{\theta,\lambda,\ell,\tau_*} N_\rho \le \frac{\varepsilon}{2}(E_\rho+D_\rho)+C Q^{m_p}.
\]
Combining this with \eqref{121} and the preceding absorption proves \eqref{76}.
\end{proof}

\section{A uniform local energy bound and the Liouville theorem} \label{123}

Assume $1<p<p_S(n,\alpha)$, so $\delta>0$. Fix $a<b$. Let $v$ be a nonnegative classical solution of \eqref{25} on a neighborhood of a compact cylinder containing the closed moving supports below, for $a\le t\le b$. Write $z_+=\max\{z,0\}$ and set
\begin{equation}\label{124}
\begin{gathered}
    \chi=(2\delta(t-a)-|x|^2)_+^5,\qquad \varphi=10( 2\delta(t-a) - |x|^2)_+^4, \qquad w= 80 |x|^2( 2\delta(t-a)-|x|^2 )_+^3,\\
    \nabla \chi = -x\varphi,\qquad w = -r\varphi_r, \qquad \Delta\chi=w-n\varphi,\qquad \chi_t = \delta \varphi.
\end{gathered}
\end{equation}
Lemma~\ref{35} applies, and the annular localization of $w$ below has the regularity required by Proposition~\ref{75}. Here double integrals are over $\R^n\times(a,b)$.

\begin{lemma}\label{125}
    For the weights in \eqref{124},
    \begin{equation}\label{126}
    \begin{aligned}
        &\iint(b-t)^3\chi v_t^2 + \frac{3(p-1)}{2(p+1)} \iint(b-t)^2 \chi |x|^\alpha v^{p+1}\\
        & \quad = \iint(b-t)^3 w\left(v_r^2-\frac{p-1}{2(p+1)}|x|^\alpha v^{p+1}\right) + \frac{1}{2}\iint(b-t)^3(w-\delta\varphi)vv_t\\
        & \qquad \quad + \frac{1}{4}\iint\Bigl\{ (b-t)^3 [(n-2)\varphi_t+2\kappa\Delta\varphi-\Delta w]\\
        & \hspace{40mm} + 6(b-t) \chi - 3(b-t)^2 [w+(\delta-2) \varphi] \Bigr\}v^2.
    \end{aligned}
    \end{equation}
    Moreover,
    \begin{equation}\label{127}
        \iint(b-t)^3 \chi |\nabla v|^2 = \iint(b-t)^3\chi |x|^\alpha v^{p+1} + \frac{1}{2} \iint\{-3(b-t)^2 \chi + (b-t)^3[w+(\delta-n) \varphi]\}v^2.
    \end{equation}
\end{lemma}

\begin{proof}
    Put $q=2\delta(t-a)-|x|^2$. Direct differentiation gives
    \[
    \begin{aligned}
        \varphi_t&=80\delta q_+^3,\\
        \Delta\varphi&=-80nq_+^3+480|x|^2q_+^2,\\
        \Delta w&=160nq_+^3 - 480(n+4)|x|^2q_+^2 + 1920|x|^4q_+.
    \end{aligned}
    \]
    Here $\chi,\varphi,w$ are $C^4,C^3,C^2$, respectively, across $q=0$. With the notation of \eqref{124}, define
    \begin{equation}\label{128}
        I(t) = \int \chi \left( \frac{1}{2} |\nabla v|^2 -\frac{|x|^\alpha v^{p+1}}{p+1} \right) + \frac{n-2}{4} \int\varphi v^2.
    \end{equation}
    All spatial integrals have compactly supported coefficients.

    For $a<t<b$, we first show that
    \begin{equation}\label{129}
    \begin{aligned}
        I'(t)+\int\chi v_t^2 & = \int w(v_r^2-\frac{p-1}{2(p+1)} |x|^\alpha v^{p+1}) + \frac{1}{2}\int(w-\delta\varphi)vv_t \\
        & \quad + \frac{1}{4} \int \bigl[(n-2) \varphi_t + 2 \kappa \Delta\varphi-\Delta w\bigr]v^2.
    \end{aligned}
    \end{equation}
    Differentiation and the equation give
    \begin{equation}\label{130}
    \begin{aligned}
        I'+\int\chi v_t^2 & = \int\varphi v_t \left( x \cdot \nabla v + \frac{n-2}{2} v \right)\\
        & \quad + \int \chi_t \left( \frac{1}{2} |\nabla v|^2 - \frac{ |x|^\alpha v^{p+1}}{p+1} \right) + \frac{n-2}{4} \int\varphi_tv^2.
    \end{aligned}
    \end{equation}
    For $C^{2,1}$ solutions, convolve in space on a slightly larger compact cylinder: $v^\nu=\varrho_\nu*v$, where $0\le\varrho_\nu\in C_c^\infty(B_\nu)$ and $\int\varrho_\nu=1$. Then $v^\nu,\nabla v^\nu,v_t^\nu$ converge uniformly to $v,\nabla v,v_t$, and
    \[
        v_t^\nu - \Delta v^\nu = |x|^\alpha(v^\nu)^p + r_\nu,\qquad r_\nu= \varrho_\nu*(|x|^\alpha v^p)-|x|^\alpha(v^\nu)^p \to 0
    \]
    uniformly. Testing with $\chi v_t^\nu$ adds the term $\int\chi r_\nu v_t^\nu$, whose time integral tends to zero. Passing to the limit gives the integrated identity. Since the right-hand side is continuous, we obtain \eqref{130}.

    The equation implies
    \begin{equation}\label{131}
        \frac{1}{2}|\nabla v|^2-\frac{|x|^\alpha v^{p+1}}{p+1} = \frac{p-1}{2(p+1)}|x|^\alpha v^{p+1} - \frac{1}{2}vv_t+\frac{1}{4}\Delta(v^2).
    \end{equation}
    Since $\chi_t=\delta\varphi$ and $ \frac{n+\alpha}{p+1}-\frac{n-2}{2}=\frac{p-1}{2(p+1)}\delta $,
    \[
        \int\chi_t\left(\frac{1}{2} |\nabla v|^2-\frac{|x|^\alpha v^{p+1}}{p+1} \right) = \frac{(p-1)\delta}{2(p+1)} \int \varphi |x|^\alpha v^{p+1}-\frac{\delta}{2} \int \varphi v v_t + \frac{\delta}{4} \int \Delta \varphi \,v^2.
    \]
    Apply \eqref{36} to the first term in \eqref{130}, at the fixed time $t$. The reaction terms cancel because their coefficients sum to zero. The identity $ n- 2 + \delta = 2\kappa$ then gives \eqref{129}.

    Local regularity gives
    \[
        |I(t)|\le C_v(t-a)^{4+\frac{n}{2}},\qquad \int\chi v^2\le C_v(t-a)^{5+\frac{n}{2}}\quad(t\downarrow a).
    \]
    The solution-dependent constant $C_v$ is used only for the endpoint limits. At $b$, the time weights are zero. Therefore
    \[
        [(b-t)^3I]_a^b=0,\qquad \left[(b-t)^2\int\chi v^2\right]_a^b=0.
    \]
    Using \eqref{131} in $I$ and integrating the $vv_t$ term in time yields
    \begin{equation}\label{132}
    \begin{aligned}
        -3\int_a^b(b-t)^2I(t)\dd t &=-\frac{3(p-1)}{2(p+1)}\iint(b-t)^2\chi |x|^\alpha v^{p+1}\\
        &\quad+\frac{1}{4}\iint \{6(b-t)\chi-3(b-t)^2[w+(\delta-2)\varphi]\}v^2.
    \end{aligned}
    \end{equation}
    Indeed, $\Delta\chi+(n-2)\varphi=w-2\varphi$ and $\chi_t=\delta\varphi$. Multiplying \eqref{129} by $(b-t)^3$, integrating in time, and using \eqref{132} gives \eqref{126}. Finally, testing the equation with $(b-t)^3\chi v$ and integrating in space and time gives \eqref{127}, because $\chi_t+\Delta\chi=w+(\delta-n)\varphi$. The endpoint term $[(b-t)^3\int\chi v^2]_a^b$ is zero by the same estimates.
\end{proof}

For $0 \le \sigma \le 1$, define
\begin{equation}\label{133}
    a_\sigma = -1-\frac{2}{\delta} - \sigma, \qquad b_\sigma=2+\sigma, \qquad R_*=2+\sqrt{10\delta+4}.
\end{equation}
From now on, set
\[
    J=(a_1-1,b_1+1) = \left(-3-\frac{2}{\delta},\,4\right).
\]
Suppose $v$ is a nonnegative classical solution on $B_{R_*} \times J$ with
\begin{equation}\label{134}
    \sup_{t\in J}\int_{B_r}v(x,t)\dd x \le C_0r^{n-\kappa} \qquad (0<r\le R_*).
\end{equation}
Use the weights in \eqref{124} with $a=a_\sigma$, denoted by $\chi_\sigma,\varphi_\sigma,w_\sigma$, and put
\begin{equation}\label{135}
    \mathcal E(\sigma) = \int_{a_\sigma}^{b_\sigma} \int (b_\sigma-t)^3 \chi_\sigma (v_t^2+|\nabla v|^2+|x|^\alpha v^{p+1} ).
\end{equation}
The closed moving supports lie strictly inside the fixed cylinder, since $ [a_1,b_1]\Subset J $ and $ 2 \delta( b_1-a_1 ) = 10 \delta + 4 < R_*^2$. The energy is nondecreasing in $\sigma$, and $\mathcal E(1)<\infty$ by local regularity.

\begin{lemma}\label{136}
    With
    \[
        m= m_p\max\left\{1,\frac{\kappa-n}{2}\right\},
    \]
    there is $C=C(n,p,\alpha)$ such that
    \begin{equation}\label{137}
        \mathcal E(\sigma)\le\vartheta\mathcal E(\Sigma) + C \vartheta^{-m_p} (\Sigma-\sigma)^{-(8m_p+m)}, \quad 0\le\sigma<\Sigma\le1,\quad 0<\vartheta<1.
    \end{equation}
\end{lemma}

\begin{proof}
    At a time $a_\sigma<t<b_\sigma$, define the positive radii $r_2 = r_2(t) $ and $ r_1 = r_1(t) $ by
    \[
        r_2^2=2\delta(t-a_\sigma),\qquad r_1^2=r_2^2+\delta (\Sigma-\sigma).
    \]
    The intermediate ball contains the inner support and satisfies
    \begin{equation}\label{138}
        r_1 \ge \sqrt{ \delta (\Sigma-\sigma)}, \qquad r_1-r_2=\frac{\delta(\Sigma-\sigma)}{r_1+r_2} \ge \frac{\delta}{2R_*} (\Sigma-\sigma).
    \end{equation}
    On $B_{r_1(t)}$,
    \[
        2\delta( t-a_\Sigma) - |x|^2\ge\delta(\Sigma-\sigma),\qquad b_\Sigma-t\ge\Sigma-\sigma.
    \]
    Thus $(b_\Sigma-t)^3\chi_\Sigma \ge \delta^5 (\Sigma-\sigma)^8 $, and
    \begin{equation}\label{139}
        \int_{a_\sigma}^{b_\sigma} \int_{B_{r_1(t)}}(v_t^2 + |\nabla v|^2 + |x|^\alpha v^{p+1}) \le C( \Sigma-\sigma)^{-8}\mathcal E(\Sigma).
\end{equation}

All coefficients in \eqref{126} and \eqref{127} are uniformly bounded on the fixed cylinder. The radius bounds in \eqref{138} imply $r_1^{2(n-\kappa)}\le C(\Sigma-\sigma)^{-(\kappa-n)_+}$ and $r_1^{n-2\kappa}\le C(\Sigma-\sigma)^{-(\kappa-\frac{n}{2})_+}$. For $0<\varepsilon\le1$, \eqref{29} therefore gives
\begin{equation}\label{140}
\begin{aligned}
    \int_{B_{r_1}}v^2 & \le\varepsilon \int_{B_{r_1}} |\nabla v|^2 + C \bigl[\varepsilon^{-\frac{n}{2}}(\Sigma-\sigma)^{-(\kappa-n)_+} + (\Sigma-\sigma)^{ - (\kappa-\frac{n}{2})_+}\bigr].
\end{aligned}
\end{equation}
Young's inequality followed by \eqref{140}, with parameter $c\varepsilon^2$ for a sufficiently small fixed $c>0$, gives
\begin{equation}\label{141}
\begin{aligned}
    C\int_{B_{r_2}}|vv_t| & \le \varepsilon \int_{B_{r_1}}(v_t^2+|\nabla v|^2) + C\bigl[\varepsilon^{-(n+1)}(\Sigma-\sigma)^{-(\kappa-n)_+} + \varepsilon^{-1} (\Sigma-\sigma)^{-(\kappa-\frac{n}{2})_+}\bigr].
\end{aligned}
\end{equation}
The same estimates hold after integration in time, up to fixed changes in the constants, since $ b_\sigma-a_\sigma\le 5 + \frac{2}{\delta} $.

To estimate the pressure, choose $r_0 = c \sqrt \varepsilon$ with a sufficiently small fixed $c>0$. Since $0 \le w_\sigma\le C|x|^2$, the pressure term over $|x| \le 3 r_0$ satisfies
\begin{equation}\label{142}
    \int_{|x|\le3r_0} w_\sigma( v_r^2-\frac{p-1}{2(p+1)} |x|^\alpha v^{p+1} ) \le Cr_0^2\int_{B_{r_2}}|\nabla v|^2 \le \varepsilon \int_{B_{r_1}}|\nabla v|^2.
\end{equation}
We discarded a nonpositive reaction term. If $r_2\le 3 r_0$, this already estimates the entire pressure. Otherwise, multiply $w_\sigma$ by a smooth radial cutoff that is zero on $B_{2r_0}$ and one outside $B_{3r_0}$, obtaining $w_0$ with
\[
    0 \le w_0 \le w_\sigma, \qquad \|w_0\|_{C^2}\le C(1+r_0^{-2}) \le C \varepsilon^{-1}.
\]
The discarded part satisfies
\[
    0 \le w_\sigma - w_0 \le w_\sigma \mathbf1_{B_{3r_0}},
\]
so, after dropping its nonpositive reaction term,
\[
    \int(w_\sigma-w_0) \left(v_r^2-\frac{p-1}{2(p+1)} |x|^\alpha v^{p+1} \right) \le Cr_0^2 \int_{B_{r_2}} |\nabla v|^2.
\]
At this fixed time, apply Proposition~\ref{75} with
\[
    \rho=r_1,\qquad h=\min\left\{1,\frac{r_1-r_2}{7}\right\},
\]
and outer radius $\rho-6h$ for the annulus. Then
\[
    \rho-6h\ge r_2+h>r_2,\qquad w_0\in C_c^2(\{r_0<|x|<\rho-6h\}).
\]
The mass bound \eqref{134} holds on $J$ for this fixed radius $\rho$. The quantity in \eqref{77} satisfies
\[
    Q\le C\varepsilon^{-1}(\Sigma-\sigma)^{-\max\{1,\frac{\kappa-n}{2}\}},
\]
because $h^{-1}\le C(\Sigma-\sigma)^{-1}$, $\rho^{-1}\le C(\Sigma-\sigma)^{-\frac{1}{2}}$, and $C_0\rho^{n-\kappa}\le C(\Sigma-\sigma)^{-\frac{(\kappa-n)_+}{2}}$. This bound is uniform in $t$. Proposition~\ref{75} applies at each fixed time with the same exponent $m_p$; no time derivative of $w_0,h,\rho$ is used. It gives
\[
\begin{aligned}
    \int w_0\left(v_r^2-\frac{p-1}{2(p+1)}|x|^\alpha v^{p+1}\right)
    &\le\varepsilon\int_{B_{r_1}}(v_t^2+|\nabla v|^2+|x|^\alpha v^{p+1})\\
    &\quad+C\varepsilon^{-m_p}(\Sigma-\sigma)^{-m_p\max\{1,\frac{\kappa-n}{2}\}}.
\end{aligned}
\]
Since $m_p\ge n+1$, the definition of $m$ gives
\[
    m \ge \max \left\{(\kappa-n)_+,(\kappa-\frac{n}{2})_+, m_p \max\{ 1,\frac{\kappa-n}{2} \}\right\}.
\]
The powers of $\varepsilon^{-1}$ in \eqref{140}--\eqref{141} do not exceed $m_p$, and their powers of $(\Sigma-\sigma)^{-1}$ do not exceed $m$. Multiplying by $(b_\sigma-t)^3$ and adjusting the small parameter by a fixed factor gives
\begin{equation}\label{143}
\begin{aligned}
    \int_{a_\sigma}^{b_\sigma}\int  (b_\sigma-t)^3w_\sigma(v_r^2-\frac{p-1}{2(p+1)}|x|^\alpha v^{p+1}) & \le \varepsilon\int_{a_\sigma}^{b_\sigma}  \int_{B_{r_1(t)}} (v_t^2+|\nabla v|^2 + |x|^\alpha v^{p+1}) \\
    & \quad + C\varepsilon^{-m_p}(\Sigma-\sigma)^{-m}.
\end{aligned}
\end{equation}

We now combine the two integral identities, using $\iint$ for integration over $\R^n\times(a_\sigma,b_\sigma)$. Since
$0<b_\sigma-t\le5+\frac{2}{\delta}$, we have
\[
    \frac{3(p-1)}{2(p+1)}(b_\sigma-t)^2 \ge\frac{3(p-1)}{2(p+1)(5+\frac{2}{\delta})}(b_\sigma-t)^3.
\]
All coefficients of $v^2$ and $vv_t$ in \eqref{126} are uniformly bounded and are zero outside $B_{r_2(t)}$. It follows that
\begin{equation}\label{144}
\begin{aligned}
    & c \iint(b_\sigma-t)^3 \chi_\sigma (v_t^2 + |x|^\alpha v^{p+1})\\
    \quad \le & \iint(b_\sigma-t)^3 w_\sigma \left(v_r^2-\frac{p-1}{2(p+1)} |x|^\alpha v^{p+1} \right) + C \int_{a_\sigma}^{b_\sigma}\int_{B_{r_2(t)}}(v^2+|vv_t|).
\end{aligned}
\end{equation}
Here $c>0$ and $C$ depend only on $n,p,\alpha$. Similarly, \eqref{127} gives
\begin{equation}\label{145}
\begin{aligned}
    \iint(b_\sigma-t)^3\chi_\sigma|\nabla v|^2 &\le\iint(b_\sigma-t)^3\chi_\sigma|x|^\alpha v^{p+1}\\
     & \quad+ C \int_{a_\sigma}^{b_\sigma} \int_{B_{r_2(t)}} v^2.
\end{aligned}
\end{equation}
Insert \eqref{143} and the time-integrated estimates \eqref{140}--\eqref{141}. All these estimates use $B_{r_1(t)}$. Adjusting the small parameters by fixed factors gives
\begin{equation}\label{146}
\begin{aligned}
    \mathcal E(\sigma) & \le C \varepsilon \int_{a_\sigma}^{b_\sigma} \int_{B_{r_1(t)}} (v_t^2+|\nabla v|^2+|x|^\alpha v^{p+1})\\
    & \quad + C\varepsilon^{-m_p} (\Sigma-\sigma)^{-m}.
\end{aligned}
\end{equation}
These fixed factors depend only on the equation parameters and do not change the exponents. Use \eqref{139} and choose $\varepsilon=c\vartheta (\Sigma-\sigma)^8$ with a sufficiently small fixed $c$.
The remainder then has the bound
\[
    C\varepsilon^{-m_p}(\Sigma-\sigma)^{-m} \le C\vartheta^{-m_p}(\Sigma-\sigma)^{-(8m_p+m)}. 
\]
\end{proof}

\begin{proposition}\label{147}
Under \eqref{134},
\begin{equation}\label{148}
    \int_{-1}^{1}\int_{B_1}(v_t^2+|\nabla v|^2+|x|^\alpha v^{p+1})  \le C(n,p,\alpha).
\end{equation}
The constant is independent of $\mathcal E(1)$ and of the supremum of $v$.
\end{proposition}

\begin{proof}
Set
\[
    \sigma_j=1-2^{-j}\quad(j\ge0),\qquad \vartheta=2^{-8m_p-m-2}.
\]
Since $\sigma_{j+1}-\sigma_j=2^{-j-1}$, iterating \eqref{137} gives
\[
    \mathcal E(0) \le\vartheta^j\mathcal E(\sigma_j) + C\vartheta^{-m_p} 2^{8m_p+m} \sum_{i=0}^{j-1}(\vartheta2^{8m_p+m})^i.
\]
Here $ \vartheta2^{8m_p+m} = \frac{1}{4} $ and
$ \vartheta^j \mathcal E(\sigma_j) \le \vartheta^j \mathcal E(1) \to 0$. Thus $\mathcal E(0)\le C(n,p,\alpha)$. On $B_1\times[-1,1]$,
\[
    2 \delta(t-a_0) - |x|^2 \ge 3,\qquad b_0 - t \ge 1.
\]
Therefore
\[
    3^5 \int_{-1}^{1} \int_{B_1} (v_t^2+|\nabla v|^2+|x|^\alpha v^{p+1}) \le \mathcal E(0)\le C(n,p,\alpha).\qedhere
\]
\end{proof}

The uniform local estimate completes the scaling argument.

\begin{proof}[Proof of Theorem~\ref{8}]
    For $R>0$, set $v_R(y,s)=R^\kappa u(Ry,R^2s)$. By \eqref{24},
    \[
        \sup_s\int_{B_r}v_R(y,s)\dd y\le C_0r^{n-\kappa} \qquad(0<r\le R_*),
    \]
    with $C_0$ independent of $R$. Proposition \ref{147} therefore applies uniformly in $R$. For each fixed $R$, the outer energy is finite, and the iteration removes its contribution.

    Changing variables in \eqref{148}, with $Q_R=B_R\times(-R^2,R^2)$, gives
    \begin{equation}\label{149}
        R^\delta \iint_{Q_R}u_t^2 + R^{\delta-2}\iint_{Q_R}(|\nabla u|^2 + |x|^\alpha u^{p+1}) \le C.
    \end{equation}
    Here $\kappa(p+1)-\alpha-n-2=\delta-2$. For each fixed $r>0$ and all $R\ge r$,
    \[
        \iint_{Q_r}u_t^2\le CR^{-\delta}\longrightarrow0.
    \]
    Thus $u_t=0$, so $u$ is independent of time. This step precedes the spatial estimate because $\delta-2$ need not be positive. Consequently,
    \[
        2R^\delta\int_{B_R} (|\nabla u|^2+|x|^\alpha u^{p+1})\le C.
    \]
    For fixed $r>0$, the integral over $B_r$ is at most $\frac{C}{2}R^{-\delta}$ for every $R\ge r$, and hence is zero. Therefore $ u=0$ away from the origin and, by continuity, also at the origin.
\end{proof}

\begin{proof}[Sketch of the proof of Theorem~\ref{9}]
We follow the proof of Theorem~\ref{8} and indicate the changes needed for the coupled nonlinearity. Set
\[
    F(U)=\frac{1}{p+1}\sum_{i,j=1}^m \beta_{ij}u_i^{r+1}u_j^{r+1},\qquad p=2r+1.
\]
The symmetry and sign assumptions on the coefficients give
\[
    U_t-\Delta U=|x|^\alpha\nabla F(U),\qquad U\cdot\nabla F(U)=(p+1)F(U),\qquad c|U|^{p+1}\le F(U)\le C|U|^{p+1}.
\]
Since $(u_i)_t-\Delta u_i\ge\beta_{ii}|x|^\alpha u_i^p$,
the proof of Lemma \ref{16} applies to each component. Consequently,
\[
    \sup_t\int_{B_R}\sum_{i=1}^m u_i(x,t) \dd x \le C_0R^{n-\kappa},\qquad \kappa=\frac{\alpha+2}{p-1}.
\]
This estimate is preserved by the rescaling
$V(y,s)=R^\kappa U(Ry,R^2s)$, and Lemma \ref{27} applies componentwise.

For the local pressure estimate, we use a vector version of
Proposition \ref{75}, with
\[
    E_\rho = \int_{B_\rho} \bigl( |\nabla V|^2 + |x|^\alpha F(V) \bigr), \qquad D_\rho=\int_{B_\rho}|V_t|^2.
\]
For the annular weights used in Section~\ref{123}, this estimate reads
\[
    \int w_0\left(|V_r|^2-\frac{p-1}{2}|x|^\alpha F(V)\right) \le\varepsilon(E_\rho+D_\rho)+CQ^{M_*}.
\]
The proof of this estimate uses the superlevel sets of the components rather than those of the coupled sources. The comparison under interval dilation treats the mixed monomials in $F$ through intersections of the component superlevel sets.

The localized Pohozaev and moving-weight identities follow by summing over the components and using the homogeneity identity above. In these identities, $v^2$, $v_t^2$, $|\nabla v|^2$, and $|x|^\alpha v^{p+1}$ are replaced by $|V|^2$, $|V_t|^2$, $|\nabla V|^2$, and $(p+1)|x|^\alpha F(V)$, respectively. The preceding pressure bound and the componentwise interpolation estimates therefore yield the iteration of Section \ref{123}, with a
possibly different exponent in the polynomial remainder. Hence
\[
    \int_{-1}^{1}\int_{B_1} \bigl(|V_t|^2+|\nabla V|^2+|x|^\alpha F(V)\bigr) \dd x \dd t\le C.
\]
Returning to the original variables gives
\[
    R^\delta\iint_{Q_R}|U_t|^2 + R^{\delta-2} \iint_{Q_R} \bigl( |\nabla U|^2 + |x|^\alpha F(U) \bigr) \le C, \qquad \delta = 2 \kappa + 2 - n > 0.
\]
Letting $ R \to \infty $ first gives $U_t=0$. The time integration then contributes a factor $ 2 R^2 $, and a second limit gives $F(U)=0$ away from the origin. The lower bound for $F$ and continuity imply $ U \equiv 0 $. All constants may depend on $m$ and $(\beta_{ij})$, but not on $R$ or the particular solution.
\end{proof}

\medskip
\noindent\textbf{Acknowledgments}
Yuxia Guo was supported by the National Natural Science Foundation of China (No. 12271283) and National Key R\&D Program (2023YFA1010002).

\medskip
\noindent\textbf{Statements and Declarations}
The authors confirm that there are no relevant financial or non-financial competing interests to report.

\medskip
\noindent\textbf{Data Availability Statement}
All data generated or analyzed during this study are included in this article.

\medskip
\noindent\textbf{Use of AI.}
ChatGPT (OpenAI) was used to discuss research ideas and proof strategies, assist with checks of mathematical arguments and notation, and support language editing, LaTeX preparation, and compilation. All mathematical arguments were independently verified by the authors, who take full responsibility for the results and the final manuscript.

\end{document}